\documentclass[a4paper,10pt]{article}
\usepackage[hmarginratio=1:1, bottom=3.5cm, top=2.5cm]{geometry}
\usepackage{amssymb,amsthm,amsmath,latexsym,url,graphicx, comment}
\usepackage[english]{babel} 
\usepackage[dvipsnames]{xcolor}
\usepackage[utf8]{inputenc}
\usepackage[sc,osf]{mathpazo}
\usepackage[euler-digits]{eulervm}
\usepackage[T1]{fontenc}
\usepackage{mathtools}
\usepackage{mathtools}
\usepackage{pdfpages}
\usepackage{float}
\usepackage[colorlinks=true,pdftex,unicode=true,linktocpage,bookmarksopen,hypertexnames=false]{hyperref}
\theoremstyle{dotless}
\newtheorem{theo}{Theorem}

\newtheorem{lem}[theo]{Lemma}
\newtheorem{cor}[theo]{Corollary}
\newtheorem{rem}[theo]{Remark}
\newtheorem{ex}[theo]{Example}
\newtheorem{defn}[theo]{Definition}
\newtheorem{quest}[theo]{Question}
\counterwithin{theo}{section}

\newcommand{\add}[1]{(#1, +)}
\newcommand{\mul}[1]{(#1, \circ)}

\newcommand{\muladd}[1]{(#1, \circ, +)}

\title{On Dedekind Skew Braces}
\author{A. Caranti -- I. Del Corso -- M. Di Matteo -- M. Ferrara -- M. Trombetti}

\usepackage{tikz-cd}

\begin{document}
\date{}
\maketitle

\begin{abstract} 
\noindent Skew braces play a central role in the theory of set-theoretic non-degenerate solutions of the Yang--Baxter equation, since their algebraic properties significantly affect the behaviour of the corresponding solutions (see for example [Ballester-Bolinches et al., {\it Adv. Math.} 455 (2024), 109880]). Recently, the study of nilpotency-like conditions for the solutions of the Yang--Baxter equation has drawn attention to skew braces of abelian type in which every substructure is an ideal (so-called, {\it Dedekind} skew braces); see for example [Ballester-Bolinches et al., {\it Result Math.} 80 (2025), Article Number 21].

The aim of this paper is not only to show that the hypothesis that the skew brace is of abelian type can be dropped in essentially all the known results in this context, but also to extend this theory to skew braces whose additive or multiplicative groups are locally cyclic (and, more generally, of finite rank).

Our main results --- which are in fact much more general than stated here --- are as follows:
\begin{itemize}
    \item Every finite Dedekind skew brace is centrally nilpotent.
    \item Every hypermultipermutational Dedekind skew brace with torsion-free additive group is trivial.
    \item We characterize skew braces whose additive or multiplicative group is locally cyclic.
    \item If a set-theoretic non-degenerate solution of the Yang--Baxter equation whose structure skew brace is Dedekind and fixes the diagonal elements, then such a solution must be the twist solution.
\end{itemize}
\end{abstract}

\medskip\medskip

central nilpotency; locally cyclic skew brace

\bigskip
\noindent {\bf 2020 Mathematics Subject Classification}\\
{\it Primary}: 16T25, 16N40, 81R50\\
{\it Secondary}: 20F16, 20F18, 20K15 

\section{Introduction}



\medskip

The Yang--Baxter equation (YBE, for short) is a consistency equation which was independently obtained by the physicists Yang \cite{Yang} and Baxter \cite{Baxter} in the field of quantum statistical mechanics. It has many relevant interpretations in the realm of mathematical physics, and besides that, it plays a key role in the foundation of quantum groups. Moreover, it provides a multidisciplinary link between a wide variety of areas such as Hopf algebras, knot theory, and braid theory, among others. In order to tackle the challenge of identifying all the so-called set-theoretic, non-degenerate solutions (commonly referred to as ‘‘solutions’’) to the YBE, the concept of skew brace has recently been introduced in \cite{GV} as a generalization of the left braces defined by Rump in \cite{Rump0}.

A {\it skew} ({\it left}) {\it brace} is a triple $(B, +, \circ)$ where $(B, +)$ and $(B, \circ)$ are (not necessarily abelian) groups such that the {\it skew \textnormal(left\textnormal) distributivity} holds: $$a \circ (b+c)=a\circ b-a+a \circ c,$$ for all $a, b, c \in B$. Thus, skew braces $(B,+,\circ)$ of {\it abelian type} (\hbox{i.e.}, such that $(B,+)$ is abelian) are the braces introduced by Rump. Clearly, the algebraic structure of a group provides the easiest examples of skew braces: If $(G,\cdot)$ is a group, then $(G,\cdot,\cdot)$ and $(G,\cdot,\cdot^{op})$ (where $g\cdot^{\operatorname{op}}h=h\cdot g$ for $h,g\in G$) are skew braces called respectively \textit{trivial} and \textit{almost trivial} (in the following, we will usually use the word ‘‘group’’ as a synonym for ‘‘trivial skew brace’’).

It turns out that to any solution of the YBE, one can associate a skew brace (sometimes called the {\it structure skew brace}) whose structure provides significant information about the solution. For example, if the structure skew brace satisfies a certain nilpotency condition, then the solution is a well-behaved solution (in the sense that it is multipermutational). As shown by many papers (see \cite{centnilpo},\cite{supersol},\cite{ballester2024soluble},\cite{edinburgh},\cite{cedo1},\cite{nilpotent},\cite{cindy},\,\ldots), the study of nilpotency and solubility concepts of skew braces is a key area of research.

Now, under suitable nilpotency conditions, a skew brace  usually has many ideals (that is, substructures that allow one to take quotients), and this fact was exploited in \cite{dedekindbraces} to introduce and study an essentially stronger concept of nilpotency of skew braces. In fact, \cite{dedekindbraces} focuses on the study of skew braces {\it of abelian type} in which {\it all} substructures are ideals (these have been called Dedekind braces), and the authors are able to prove for example that every finite skew brace of abelian type satisfies one of the strongest possible nilpotency conditions so far introduced for skew braces: central nilpotency. They are also able to characterize certain types of Dedekind braces. The same authors continued investigating this fruitful nilpotency condition in a subsequent paper \cite{ballesterDedekindii}, where they determined the influence of a  Dedekind structure skew brace of abelian type on the corresponding involutive solution of the YBE.

The aim of our paper is not only to show how to drop the 'abelian type' assumption from virtually all the results in the previously mentioned papers, but also to promote the study of Dedekind skew braces so as to encompass a wide variety of both finite and infinite skew braces.  In doing so, we also obtain, as a by-product, several fruitful characterizations of certain classes of skew braces. We now briefly outline our main results and their consequences. 

\begin{itemize}
    \item {\it Every finite Dedekind skew brace is centrally nilpotent \textnormal(with no uniform bound on the class\textnormal).}\\ 
    This is the main content of Section \ref{sec1vero}. Actually, Theorem \ref{ThA} and its corollaries describe the situation in a much broader range of situations (especially when one of the underlying groups is non-periodic). For example, Corollary \ref{torsionfreeincredible} shows that when the additive group is torsion-free, any hypermultipermutational Dedekind skew brace is a trivial brace (see~Co\-rol\-la\-ry~\ref{solution1} concerning the solutions of the YBE), while Corollary \ref{classe2} gives a good description of Dedekind skew braces of multipermutational level $2$. Examples and open questions are given at the end of Section \ref{sec1vero}.

    \item {\it A skew brace of prime power order is Dedekind provided that either the additive or the multiplicative group is cyclic.}\\
    This is the content of Theorem \ref{thB}, and implicitly provides many examples of non-trivial~De\-de\-kind skew braces. The subsequent results (Corollary \ref{corfinitecycliccase}, Theorems \ref{cyclicideals} and~\ref{cyclicsupersoluble}) give some more information in the composite case, also showing how Dedekind-ness can be used to simplify proofs of certain known results.

    \item {\it We provide necessary and sufficient conditions for a skew brace to have either an additive or multiplicative torsion-free locally cyclic underlying group.}\\
    This is achieved in Theorems \ref{charQ+} and \ref{charQcirc}. A consequence of these descriptions is that in the Dedekind case such skew braces must be trivial (see Corollary \ref{thinfinitecyclic}). Two further consequences (see Theorems \ref{thpolycyclic} and \ref{divisibleideals}) of this characterization deal with skew braces whose multiplicative or additive groups satisfy certain finite rank conditions.  
    
\item {\it Let $(X,r)$ be a solution of the YBE such that $G(X,r)$ is a Dedekind skew brace. If $r(x,x)=(x,x)$ for every $x\in X$, then $(X,r)$ is a twist solution}.\\
This is Corollary \ref{corsolution}, which depends on Theorem \ref{theoxastx}. A stronger and solution-related form of Dedekind-ness is examined in Theorems \ref{generalizeadolf} and \ref{genadolfo2}.
    
\end{itemize}

\medskip


We have kept the introduction free of excessive technical details in order to focus on the meaning and motivation behind our approach. Readers are referred to the next section for the background on notation, definitions, and structural results.

\section{Preliminaries}\label{preliminarysec}

Since skew braces are a very recent tool for studying the Yang–Baxter equation, the terminology and notation have not yet been standardized. There are two main (completely equivalent) approaches here: considering {\it right} or {\it left} skew braces, and this typically depends on the background of the researchers. As the reader may have deduced from the introduction, this paper will only deal with {\it left} skew braces (which will be simply called {\it skew braces} in the following). But even within the realms of right and left skew braces, there exist subcultures, with different communities favoring distinct conventions and notations. For example, sometimes the first operation of a skew brace is denoted by a plus-symbol, other times by a central dot. The purpose of this section is thus to clarify our terminology for the reader, even for those already acquainted with skew braces. We will also establish some general and preliminary results, which will be frequently used throughout the remainder of the paper.

Let $(B,+,\circ)$ be a skew brace. The groups $(B,+)$ and $(B,\circ)$ are respectively known as the {\it additive} and {\it multiplicative} groups of the skew brace. The identity element of $(B,+)$, denoted by~$0$, and the identity element of $(B,\circ)$, denoted by $1$, are the same element, and we generally use~$0$ to denote this common identity element, even in a multiplicative context. Let $b\in B$. If $m$ is any positive integer, we put $$mb=\underbrace{b+\ldots+b}_{\textnormal{$m$ times}}\quad\textnormal{and}\quad b^m=\underbrace{b\circ\ldots\circ b}_{\textnormal{$m$ times}}.$$ If $m$ is any negative integer, we put $mb=-\big((-m)b\big)$ and $b^{m}=\big(b^{-m}\big)^{-1}$. In particular, $-b$ and $b^{-1}$ are the inverses of $b$ in $(B,+)$ and $(B,\circ)$, respectively. Clearly, we also put \hbox{$0b=0=1=b^0$.} Moreover, for any integer $m$, we define $(mB,+)$ as the additive subgroup of~$B$ generated by the elements of the form $mc$, where $c\in B$, so  \hbox{$(mB,+)\!=\!\langle mc: c\in B\rangle_+$.} More generally, for a subset $X$ of $B$, we denote by $\langle X\rangle_+$ the subgroup generated by $X$ in the group $(B,+)$. Similar definitions go for~$(B^m,\circ)$ and $\langle X\rangle_\circ$. Note that in general, if we are considering a standard term from the world of group theory with an added plus or circle subscript, we are referring to the term in the context of the additive or multiplicative group of the skew brace.

It turns out that the map
$$\lambda : a \in (B, \circ) \mapsto \big(\lambda_a:b\in B\mapsto -a+a \circ b\in B\big) \in \operatorname{Aut}(B, +),$$ is a group homomorphism (in the context of right skew braces the $\lambda$-map is usually denoted by $\gamma$).  The following formulas are easily verified
$$a+b=a \circ \lambda_a^{-1}(b), \quad a \circ b=a+\lambda_a(b), \quad -a=\lambda_a(a^{-1})$$ and will frequently be used without any further notice. Moreover, an important role is played by the so-called \textit{star operation}:
$$a \ast b=\lambda_a(b) - b=-a + a \circ b-b$$ 
for all $a,b \in B$, which, roughly speaking, "measures" how different the two operations on $B$ are. The following conditions tell us how the star operation interacts with the operations of the skew brace and they are very useful in calculations:
$$a \ast (b+c) = a \ast b +b+a \ast c -b\quad\textnormal{and}\quad(a \circ b) \ast c=a \ast (b \ast c)+b \ast c+ a \ast c$$
for all $a, b, c \in B$. If $X$ and $Y$ are subsets of $B$, then we put $X\ast Y=\langle x\ast y\,:\, x\in X,\, y\in Y\rangle_+$. The $\lambda$-map allows us to consider the following natural semidirect product: $$G=(B,+)\rtimes_\lambda(B,\circ),$$ where the operation is defined as follows: $$(a,b)\cdot(c,d)=(a+\lambda_b(c), b \circ d)$$ for all $a, b, c, d \in B$. An easy computation shows that the star operation corresponds to a commutator: 
$$[(0, a),(b, 0)] = (a \ast b, 0)$$ 
for all $a, b \in B$) (note that our convention for commutators is $[x, y]=xyx^{-1}y^{-1}$).

A very prominent part in the study of a skew brace relies on the study of the sub-structures. The main ones we are going to use are the following.

\begin{itemize}
    \item \textit{sub-skew brace}: a subgroup of both $(B,+)$ and $(B,\circ)$. The sub-skew brace generated by a subset $X$ of $B$ (that is, the smallest sub-skew brace containing $X$ with respect to the inclusion) of $B$ is denoted by $\langle X\rangle$. The property of being a sub-skew brace for a subset $F$ will be emphasized by writing $F\leq B$. The kernel $\operatorname{Ker}(\lambda)$ of the $\lambda$-map is a very relevant example of a sub-skew brace (which is obviously also normal in the multiplicative group).
    \item \textit{left-ideal}: a subgroup of $(B,+)$ which is {\it $\lambda$-invariant} (namely, invariant under the automorphism~$\lambda_b$ for every $b\in B$). A left-ideal is automatically a sub-skew brace.
    \item \textit{strong left-ideal}: a left-ideal that is normal in $(B,+)$. It is easy to see that every characteristic subgroup of $(B,+)$ is a strong left-ideal of $B$.
    \item \textit{ideal}: a strong left-ideal that is normal in $(B,\circ)$. If $I$ is an ideal of $B$, then we usually write also $I\trianglelefteq B$. Additive and multiplicative cosets coincide for ideals, and this allows for considering the quotient skew brace with the obvious inherited operations. The ideal generated by a subset $X$ of $B$ (that is, the smallest ideal containing $X$ with respect to the inclusion) is denoted by $\langle X\rangle^B$. Ideals can equivalently be thought of as kernels of skew brace {\it homomorphisms} (that is, maps preserving both operations).
\end{itemize}

There are many ways to produce (left-)ideals starting from known (left-)ideals. For example, if $L$ is a left-ideal and $I$ is an ideal, then $I\ast L$ is a left-ideal of $B$, while $I\ast B$ is always an ideal of~$B$. Thus, for example, $B\ast B$ is always an ideal of~$B$. The relevance of this ideal can for instance be seen from the fact that $B$ is trivial as a skew brace if and only if $B\ast B=\{0\}$. Two further very relevant ideals are the {\it socle} $\operatorname{Soc}(B)$ and the {\it centre} $Z(B)$ of $B$. The former is defined as the intersection of $\operatorname{Ker}(\lambda)$ with the additive centre $Z(B,+)$, while the latter is the intersection of the socle with the multiplicative centre $Z(B,\circ)$. These ideals are very relevant in the context of nilpotency conditions for a skew brace.

The following results are probably folklore and allow us to easily decide on many occasions if a subset is a sub-skew brace or an ideal. First, recall that a subset $X$ of a group $(G, \cdot)$ is called a \textit{normal} if \hbox{$g^{-1}\cdot X\cdot g \subseteq X$,} for all $g \in G$.

\begin{lem}\label{bastb}
Let $(B,+,\circ)$ be a skew brace. If $b\in B$ with $b\ast b=0$, then the sub-skew brace $\langle b\rangle$ generated by $b$ is trivial, so in particular $\langle b\rangle$ coincides with both the additive and multiplicative subgroups~$\langle b\rangle_+$ and $\langle b\rangle_\circ$ generated by $b$.
\end{lem}
\begin{proof}
First, $$0=b\ast (-b+b)=b\ast (-b)-b+b\ast b+b=b\ast (-b).$$ Thus, $b\circ (-b)=b-b=0$ and so $-b=b^{-1}$. Moreover, $$0=(b^{-1}\circ b)\ast b=b^{-1}\ast (b\ast b)+b^{-1}\ast b+b\ast b=b^{-1}\ast b$$ and then $$0=b^{-1}\ast (-b+b)=b^{-1}\ast (-b)-b+b^{-1}\ast b+b=(-b)\ast (-b).$$ Now, suppose $nb=b^n$ for some integer $n$. We easily see that 
$(\pm b)\ast (\pm nb)=0$, so $$\pm(n+1)b=\pm b+\pm nb=(\pm b)\circ (\pm nb)=b^{\pm}\circ b^{\pm n}=b^{\pm(n+1)}.$$ Therefore $\langle b\rangle=\langle b\rangle_+=\langle b\rangle_\circ$ is a trivial brace.
\end{proof}

\begin{lem}\label{lemuseful}
Let $(B,+,\circ)$ be a skew brace and let $X$ be an additive or multiplicative subgroup of~$B$. If any three of the following four properties hold:
\begin{itemize}
    \item[\textnormal{(1)}] $X$ is a normal subset of~$(B,+)$;
    \item[\textnormal{(2)}] $X$ is $\lambda$-invariant;
    \item[\textnormal{(3)}] $X$ is a normal subset of $(B,\circ)$;
    \item[\textnormal{(4)}] $X\ast B\subseteq X$;
\end{itemize} 
\noindent then $X$ is an ideal.
\end{lem}
\begin{proof}
Note that if $X$ satisfies (2) and is a multiplicative subgroup of $B$, then $X$ is also an additive subgroup since $b\circ X=b+X$ for every $b\in B$, and $\lambda_x(x^{-1})=-x$.

If $X$ satisfies (1)--(3), then $X$ is an ideal by definition. Furthermore, if $X$ satisfies (1), (2) and~(4), then $X$ is a left-ideal and  $$b^{-1}\circ x\circ b=b^{-1}\circ(x'+b)=b^{-1}+\lambda_{b^{-1}}(x')-b^{-1}\in X$$ for every $x\in X$, $b\in B$, and for some $x'\in X$, so again $X$ is an ideal.

Suppose $X$ satisfies (1), (3) and (4). Let $b\in B$ and $x\in X$. 
Since $(X,\circ)$ is normal in~$(B,\circ)$, $x'=b\circ x\circ b^{-1}\in X$. By (4), $x'\circ b=x''+b$ for some $x'' \in X$, and hence $$\lambda_b(x)=-b+b\circ x=-b+x'\circ b=-b+x''+b\in X$$ by (1). Thus $X$ satisfies (2) and is consequently an ideal.

Suppose finally that $X$ satisfies (2)--(4). Let $b\in B$ and $x\in X$. By (2), there is $x'\in X$ such that $\lambda_b(x')=x$, while by (3) there is $x''\in X$ such that $b^{-1}\circ x''\circ b=x'$. Thus, $$-x''+b+x-b=-x''+b+\lambda_b(x')-b=-x''+x''\circ b-b=x''\ast b\in X$$ by (3), and hence $b+x-b\in X$, which means that $X$ is additively normal, and so even an ideal of $B$.
\end{proof}

\medskip

In \cite{centnilpo}, the concept of index for a sub-skew brace has been defined and has been used to derive a lot of interesting nilpotency-related properties. Let $X$ be a sub-skew brace of $B$. Clearly, if $B$ is finite, then additive and multiplicative indices coincide by order considerations and~La\-grange's theorem. 
In the infinite case, this no longer seems to hold, although no examples of this have been found so far. In any case, if the additive and multiplicative indices of~$X$ in $B$ are finite and coincide, then we say that $X$ has finite index, and we denote the common index by~$|B:X|$. If they are both infinite, we say that $X$ has infinite index. In any other case, the index is not defined. In \cite{centnilpo}, it was shown that the index of a sub-skew brace is defined in many cases. We are now going to improve one of the results from that work, and we will use it in the next section. First, we need to recall some terminology.

\begin{defn}\label{defiindexpres}
{\rm A skew brace $(B,+,\circ)$ is said to be {\it index-preserving} if the following conditions are equivalent for every sub-skew brace $C$ of $B$:
\begin{itemize}
    \item[\textnormal{(1)}] $|(B, +) : (C, +)|$ is finite.
    \item[\textnormal{(2)}] $|(B, \circ) : (C, \circ)|$ is finite.
    \item[\textnormal{(3)}] $C$ has finite index in $B$.
    \item[\textnormal{(4)}] $|B:C|$ exists.
\end{itemize}

\noindent Note that the implications (4)$\implies$(3)$\implies$(1),(2) trivially hold.}
\end{defn}

\begin{defn}
{\rm Let $\mathfrak X$ be a property pertaining to skew braces (for example the property of being finite or that of being abelian). We say that:
\begin{itemize}
    \item $\mathfrak X$ is {\it local} if a skew brace $C$ has property $\mathfrak X$ provided that each finitely generated sub-skew brace is contained in a sub-skew brace with the pro\-per\-ty~$\mathfrak X$. Clearly, if the property~$\mathfrak X$ is inherited by the sub-skew braces, then $\mathfrak X$ is local if and only if the finitely generated sub-skew braces have the property $\mathfrak X$. Thus, for example, the property of being abelian is local, while that of being finite is not. In case the property is not local, it makes sense to consider the class of {\it locally $\mathfrak X$ skew braces}, that is, the class of all skew braces in which the finitely generated sub-skew braces are contained in $\mathfrak X$ sub-skew braces; thus, for example a skew brace is {\it locally finite} if every finitely generated sub-skew brace is finite. 
    \item $\mathfrak X$ is closed with respect to forming ascending series if $C$ has property $\mathfrak X$ provided that $C$ has an ascending series of sub-skew braces $$\{0\}=I_0\trianglelefteq I_1\trianglelefteq\ldots\; I_\alpha\trianglelefteq I_{\alpha+1}\trianglelefteq\ldots\; I_\mu=C$$ (here $\lambda$ is an ordinal number) whose factors $I_{\beta+1}/I_\beta$, $\beta<\mu$, have the property $\mathfrak X$.
\end{itemize}
 }
\end{defn}

\begin{lem}\label{analogcentralnilp}
The property of being index-preserving is local and closed with respect to forming ascending series.
\end{lem}
\begin{proof}
Let $(B,+,\circ)$ be a skew brace, and $C$ a sub-skew brace. Assume first that $B$ is locally index-preserving and satisfies for example (1) of Definition \ref{defiindexpres}. Set $n=|(B, +) : (C, +)|$. Suppose first that the index of $(C,\circ)$ in $(B,\circ)$ is strictly greater than $n$, and let $b_1,\ldots,b_{n+1}$ be distinct representatives for the left cosets of $(C,\circ)$ in $(B,\circ)$. Then there exists an index-preserving sub-skew brace $D$ of $B$ containing the $b_i$'s. This implies that $$n<|(D,\circ):(C\cap D,\circ)|=|(D,+):(C\cap D,+)|\leq n,$$ a contradiction. Therefore the index of $(C,\circ)$ in $(B,\circ)$ is at most $n$. A similar argument shows now that the two indices coincide. This proves that the property of being index-preserving is local.

Suppose now that $B$ has an index-preserving ideal $I$ such that $B/I$ is index-preserving as well. Then $|(I+C,+):(C,+)|=|(I,+):(C\cap I,+)|=
|(I,\circ):(C\cap I,\circ)|=|(I\circ C,\circ):(C,\circ)|$
and
$|(B,+):(I+C,+)|=|(B,\circ):(I+C,\circ)|$. Therefore $|(B,+):(C,+)|=|(B,\circ):(C,\circ)|$ and $B$ is index-preserving. Finally, assume that $B$ has an ascending series $$\{0\}=I_0\trianglelefteq I_1\trianglelefteq\ldots\; I_\alpha\trianglelefteq I_{\alpha+1}\trianglelefteq\ldots\; I_\mu=B$$ whose factors are index-preserving. We use induction on $\mu$ to prove that $B$ is index-preserving. If $\mu$ is a limit ordinal, then $B$ is the union of index-preserving sub-skew braces, so $B$ is locally index-preserving, and consequently $B$ is index-preserving. If $\mu$ is successor ordinal, then $I_{\mu-1}$ and $B/I_{\mu-1}$ are index-preserving, so $B$ is index-preserving by what we have shown above. The result is proved.~\end{proof}

\medskip

In the context of skew braces, there exist various notions of nilpotency and solubility. These concepts were developed not only to capture the intrinsic algebraic structure of skew braces but also because they provide valuable insights into certain classes of set-theoretic solutions of the Yang–Baxter equation. Below, we briefly outline some of the most significant of these notions. For a more comprehensive account, we refer the interested reader to the following works: \cite{centnilpo}, \cite{supersol}, \cite{ballester2024soluble}, \cite{edinburgh},  \cite{stefanelli}, \cite{nilpotent}, \cite{colazzo}, \cite{JVV}.

Let $(B,+,\circ)$ be a skew brace. If $B$ is a trivial brace and $(B,+)$ is abelian, then we also say that $B$ is {\it abelian}; this is obviously equivalent to requiring that $B=Z(B)=\operatorname{Soc}(B)$. The skew brace~$B$ is said to be {\it soluble} if it has a finite series of ideals $$\{0\}=B_0\leq B_1\leq\ldots\leq B_n=B$$ with abelian factors. Now, we define the {\it upper central series} and {\it upper socle series} for $B$. First, we put $Z_0(B)=\{0\}=\operatorname{Soc}_0(B)$. Next, if $\alpha$ is any ordinal number and we have already defined the terms $Z_\alpha(B)$ and $\operatorname{Soc}_\alpha(B)$, then we define the $\alpha+1$-th terms of the series as the only ideals of~$B$ such that $$Z_{\alpha+1}(B)/Z_\alpha(B)=Z\big(B/Z_\alpha(B)\big)\quad\textnormal{and}\quad \operatorname{Soc}_{\alpha+1}(B)/\operatorname{Soc}_\alpha(B)=Z\big(B/\operatorname{Soc}_\alpha(B)\big).$$ If $\alpha$ is a limit ordinal and $Z_\beta(B)$ and $\operatorname{Soc}_\beta(B)$ have already been defined for $\beta<\alpha$, then we set $$Z_\alpha(B)=\bigcup_{\beta<\alpha}Z_\beta(B)\quad\textnormal{and}\quad\operatorname{Soc}_\alpha(B)=\bigcup_{\beta<\alpha}\operatorname{Soc}_\beta(B).$$ By construction, all terms of the upper central series and the upper socle series are ideals; clearly $Z_1(B)=Z(B)$ and $\operatorname{Soc}_1(B)=\operatorname{Soc}(B)$. The last term of the upper central series (resp. the upper socle series) is the {\it hypercentre} (resp. the {\it hypersocle}). In terms of these series, we can define some nilpotency concepts on a skew brace $B$. We say that
\begin{itemize}
    \item $B$ is \textit{centrally nilpotent} if there is a non-negative integer $n$ such that $B=Z_{n}(B)$. The smallest integer with this property is called the {\it central nilpotency class} of $B$.
    \item $B$ is {\it locally centrally-nilpotent} if every finitely generated sub-skew brace is centrally nilpotent. Locally centrally-nilpotent skew braces have a good torsion theory in the sense that an element of $B$ is additively periodic if and only if it is multiplicatively periodic (see \cite{centnilpo},~The\-o\-rem 4.12). In addition, the set of all periodic elements of $(B,+)$ forms an ideal of $B$.
    \item $B$ is \textit{hypercentral} if there is an ordinal $\alpha$ such that $B=Z_{\alpha}(B)$. This is equivalent to requiring that every non-zero homomorphic image of $B$ has a non-zero centre. It has been proved in \cite{centnilpo} that if $Z(B)$ has a torsion-free additive group, then every factor of the upper central series is a torsion-free group, so when $B$ is hypercentral, the whole additive group $(B,+)$ is torsion-free.
    \item $B$ is {\it residually centrally-nilpotent} if for every $b\in B$ there is an ideal $I$ such that  $b\not\in I$ and~$B/I$ is centrally nilpotent.
    \item $B$ has \textit{finite multipermutational level} if there is a non-negative integer $n$ with $B=\operatorname{Soc}_{n}(B)$. The smallest integer with this property is called the {\it multipermutational level of $B$.}
    \item $B$ is \textit{hypermultipermutational} if there is an ordinal number $\alpha$ such that $B=\operatorname{Soc}_{\alpha}(B)$. The consideration of any non-trivial skew brace with an infinite cyclic additive group shows that the factors of the upper socle series need not have a torsion-free additive group even if the socle has a torsion-free additive group (see also the end of this section).
    
\end{itemize}

Clearly, central nilpotency implies solubility, hypercentrality and residual central-nilpotency, while (as shown in \cite{centnilpo}) hypercentrality implies local central-nilpotency. Moreover, being of finite multipermutational level implies being soluble and hypermultipermutational. Furthermore, hypercentrality implies being hypermultipermutational, while central nilpotency implies being of finite multipermutational level (as \hbox{$Z_\alpha(B)\leq\operatorname{Soc}_\alpha(B)$} for every ordinal number~$\alpha$). Two further nilpotency concepts for skew braces that we will briefly encounter in our discussion are those of left and right nilpotency: The skew brace $B$ is {\it left} (resp. {\it right}) {\it nilpotent} if there is an integer $n$ such that 
$$\underbrace{B\ast (B\ast(\ldots \ast(B\ast B))\ldots)}_{\textnormal{$n$ times}}=\{0\}\quad\Big(\textnormal{resp. }\underbrace{(\ldots((B\ast B)\ast \ldots )\ast B)\ast B}_{\textnormal{$n$ times}}=\{0\}\Big).$$ Obviously, central nilpotency implies both properties, and  on some occasions --- for example when the skew brace is of {\it nilpotent type} (that is, its additive group is nilpotent) --- the two together imply central nilpotency (see \cite{JVV}, Corollary 2.15). Moreover, being of finite multipermutational level implies being right nilpotent. Finally, we note that the ‘‘descending’’ analogue of hypercentrality has also been considered in \cite{centnilpo} but it seems that not much can be said in the Dedekind context.

\medskip

In our work we will also frequently need the following concepts that have been introduced and studied in \cite{supersol} and \cite{poly}, and that generalize some of the previous nilpotency properties.

\begin{itemize}
    \item $B$ is {\it almost polycyclic} if there is a finite series of ideals $$\{0\}=B_0\leq B_1\leq\ldots\leq B_n=B$$ such that, for all $0\leq i<n$, either $B_{i+1}/B_i\leq\operatorname{Soc}(B/B_i)$ and $B_{i+1}/B_i$ is finitely generated, or $B_{i+1}/B_i$ is finite. This is equivalent to requiring that there is a finite series in which the first factor on the right is finite and all the others are infinite cyclic groups contained in the socle of the quotient skew brace. 
    \item $B$ is {\it supersoluble} if there is a finite series of ideals $$\{0\}=B_0\leq B_1\leq\ldots\leq B_n=B$$ such that, for all $0\leq i<n$, either $B_{i+1}/B_i\leq\operatorname{Soc}(B/B_i)$ and $B_{i+1}/B_i$ is infinite cyclic, or~$B_{i+1}/B_i$ has prime order. It turns out that if $p$ is the largest prime dividing the order of a finite supersoluble skew brace, then the additive Sylow $p$-subgroup is an ideal (see~The\-o\-rem~3.11 of \cite{supersol} for the more general result).
    \item $B$ is {\it hypercyclic} if there is an ascending series of ideals whose factors satisfy the same conditions as in the supersoluble case. 
\end{itemize}

Obviously, every supersoluble skew brace is hypercyclic and almost polycyclic; moreover, every hypercentral (resp. finitely generated and centrally nilpotent) skew brace is hypercyclic (resp. supersoluble). It has been shown in \cite{supersol} that every finite skew brace of square-free order is supersoluble (more in general this is true for every skew brace whose additive and multiplicative groups have cyclic~Sy\-low subgroups). If we want non-trivial infinite examples, we may start looking at non-trivial skew braces whose additive or multiplicative groups are infinite cyclic. As shown in~\cite{rump2007classification} and \cite{stefanellotrap}, these skew braces exist and the other group is always infinite dihedral. If the additive group is infinite cyclic, then the socle has index $2$, so this skew brace is in fact supersoluble. If the multiplicative group is infinite cyclic, then the socle is zero, so this is not even almost polycyclic, although it has an ideal which is an infinite cyclic of index $2$. The importance of these additional concepts stems from the fact that their behaviour can be almost totally controlled by the finite homomorphic images. Thus, for example, an almost polycyclic skew brace whose finite quotients are centrally nilpotent is itself centrally nilpotent.

\subsection*{Skew braces and set-theoretic solutions of the YBE}

A {\it set-theoretic solution} (we will henceforth omit the adjective set-theoretic) of the Yang--Baxter equation is a pair $(X, r)$, where $X$ is a set and $$r : (x,y)\in X \times X \mapsto r(x,y):=(\lambda_x(y), \rho_y(x))\in X \times X$$ is a map such that $r_{12} r_{23} r_{12} = r_{23} r_{12} r_{23}$  where $r_{12} = (r \times \mathrm{id}_X)$ and $r_{23} = (\mathrm{id}_X \times r)$. The solution is said to be {\it involutive} if $r^2 = \mathrm{id}_{X \times X}$, and {\it non-degenerate} if $\lambda_x$ and $\rho_x$ are injective for every $x \in X$. 

The simplest example of a solution on the set $X$ is the so-called \emph{twist solution}: $r(x, y) = (y, x)$ for every $x \in X$, i.e.\ $\lambda_x = \rho_x = \mathrm{id}_X$ for every $x \in X$. Moreover, once we have been given a solution $(X,r)$, we can construct other solutions that will be somewhat useful in understanding the structure of the original solution. The \emph{retraction solution} $\mathrm{Ret}(X, r) = (X / {\sim}_r, \bar{r})$ is defined by means of the equivalence relation ${\sim}_r$ in $X$: $x \sim y$ if $\lambda_x = \lambda_y$ and $\rho_x = \rho_y$. Then,
\[
\bar{r}([x], [y]) = ([\lambda_x(y)], [\rho_y(x)]), \quad \text{for all } [x], [y] \in X / {\sim}_r.
\] The {\it retraction series} is defined recursively as
\[
\mathrm{Ret}^0(X, r) = (X, r), \qquad
\mathrm{Ret}^{\alpha+1}(X, r) = \mathrm{Ret}(\mathrm{Ret}^\alpha(X, r)), \quad \text{for every ordinal } \alpha,
\]
\[
\mathrm{Ret}^\lambda(X, r) = (X / {\sim}_\lambda, \bar{r}_\lambda), \quad \text{for every limit ordinal } \lambda,
\]
where ${\sim}_\lambda := \bigcup_{\mu < \lambda} {\sim}_\mu$ and ${\sim}_\mu$ is the equivalence relation in $\mathrm{Ret}^\mu(X, r)$ for every $\mu < \lambda$.

A solution $(X, r)$ is \emph{multipermutational of level $m$}, if $m$ is the smallest non-negative integer such that $\mathrm{Ret}^m(X, r)$ has cardinality $1$. For instance, twist solutions are multipermutational solutions of level 1. We say that $(X, r)$ is \emph{hypermultipermutational} if for some ordinal number $\alpha$, $\mathrm{Ret}^\alpha(X, r)$ has cardinality $1$. In case that $X$ is a finite set, it is clear that a solution is hypermultipermutational if and only if it is multipermutational of finite level.

For every (finite) skew brace $(B,+,\circ)$, one can naturally associate to it a (finite) solution~$(B,r_B)$ of the YBE defined by
\begin{equation}\label{rB} r_B: (a,b) \in B\times B \longmapsto (\lambda_a(b),\lambda_a(b)^{-1}\circ a \circ b)\in B\times B,\end{equation}
which is involutive if and only if $B$ is a brace (see \cite{GV}, \cite{Rump0}, and also \cite{ChildsYBE}). We shall refer to these solutions that arise from a skew brace as {\it skew-brace-solutions} ({\it sb-solution} for short). Conversely, for every (finite) solution $(X,r)$ of the~YBE, one can associate to it the (not necessarily finite) {\it structure group} of $(X,r)$, defined by the presentation
\[G(X,r)=\langle x\in X\,:\, x \circ y=u\circ v\textnormal{ for }r(x,y)=(u,v)\rangle,\]
on which one can define a group operation $+$ such that $(G(X,r),+,\circ)$ is a skew brace satisfying a certain universal property (see \cite[Theorem 3.5]{SV}). Note that if $(X, r)$ is a twist solution, then~$G(X, r)$ is an abelian brace isomorphic to $\mathbb{Z}^{(X)}$, the free abelian group on $X$. Moreover, the following diagram commutes
\[
\begin{tikzcd}\label{diagram3}
X \times X \arrow{r}{r} \arrow[swap]{d}{\iota_X \times \iota_X} & X \times X \arrow{d}{\iota_X \times \iota_X} \\
G \times G \arrow{r}{r_G} & G \times G
\end{tikzcd}
\tag{3}
\]

\noindent where $G := G(X, r)$, $(G, r_G)$ is the solution associated with $G$, and $\iota_X : X \to G(X, r)$ is the canonical inclusion. This correspondence paves the way to characterise the multipermutational character of a solution in terms of nilpotency of left braces. It follows that $(X, r)$ is hypermultipermutational (resp. has finite multipermutational level) if and only if $G(X, r)$ is hypermultipermutational (resp. has finite multipermutational level); see \cite{cedo1}, \cite{nilpotent}.


\subsection*{Some group-theoretic concepts and terminology}

The infinite cyclic group will be denoted by $(\mathbb Z,+)$. Also, if $n$ is any positive integer, then the cyclic group of order $n$ will be denoted by $(\mathbb Z_n,+)$. Let $A$ be an abelian group. Then $A$ is {\it locally cyclic} if every finitely generated subgroup of $A$ is cyclic; in this case, $A$ embeds either in the additive group of rational numbers $(\mathbb Q,+)$ or in its factor group $(\mathbb Q/\mathbb Z,+)$. In the former case,~$A$ is torsion-free (and may always be regarded as containing the number $1$), while in the latter, $A$ is periodic.  If~$p$ is any prime number, then $A$ is~{\it $p$-divisible} if for each element $a\in A$ there is $b\in A$ with $b^p=a$. We say that 
$A$ is {\it divisible} if it is~\hbox{$p$-di}\-vi\-sible for every prime~$p$. The only divisible subgroups of $(\mathbb Q,+)$ are $\{0\}$ and $\mathbb Q$ itself.

The {\it $0$-rank} of $A$ is the cardinality of a maximal independent set consisting of elements of $A$ with infinite order. If $p$ is a prime, the $p$-rank of $A$ is the cardinality of a maximal independent set of non-trivial elements of $A$ each of which has order a power of $p$. Recall also that $A$ has {\it finite abelian subgroup rank} if it has finite $0$-rank and finite $p$-rank for every prime $p$.

Let $G$ be a group having a finite series whose factors are either infinite cyclic or periodic. The number of infinite cyclic factors in such a series is an invariant of $G$ which is called the~{\it Hirsch length}. This invariant will be frequently used to make induction arguments. For example, the {\it infinite dihedral group}, \hbox{i.e.} the semidirect product $\mathbb Z\rtimes \mathbb Z_2$ where $\mathbb Z_2$ inverts~$\mathbb Z$, has~Hirsch length $1$. Clearly, in the abelian case, the Hirsch length corresponds to the $0$-rank.

Moreover, if $G$ is any group, then $\pi(G)$ denotes the set of all primes $p$ for which $G$ has an element of order $p$ --- this is similar to what one usually does with a natural number $m$, where~$\pi(m)$ is the set of all prime divisors of $m$. Clearly, if $G$ is finite, then $\pi(G)$ is precisely the set of prime divisors of the order of $G$. Also, if $G$ is any nilpotent group, then $\pi(G)=\pi(T)$, where $T$ is the {\it torsion subgroup} of $G$, that is, the subgroup made by all the periodic elements of~$G$. Recall also that if $\pi$ is any set of primes, then $\pi'$ denotes the set of primes that are not in~$\pi$, while $O_\pi(G)$ is the largest normal $\pi$-subgroup of $G$.

The group $G$ is said to be {\it minimax} if it has a finite series whose factors satisfy either the minimal or the maximal condition on subgroups. In the {\it soluble-by-finite} case (that is, when there is a normal soluble subgroup of finite index), the factors satisfying the minimal condition become {\it \v Cernikov groups} (finite extensions of a direct product of finitely many Pr\"ufer groups), while those satisfying the maximal condition become polycyclic-by-finite. Although it may appear to be very complicated, the structure of soluble-by-finite minimax groups is described by Theorem 10.33 of \cite{Ro72}: Their {\it finite residual} $R$ (that is, the intersection of all normal subgroups of finite index) is a direct product of finitely many Pr\"ufer groups, while the {\it Fitting subgroup}~$F/R$ (that is, the product of all normal nilpotent subgroups) is nilpotent, and $G/F$ is finitely generated and abelian-by-finite. Minimax groups have {\it finite rank} in the sense that the minimal number of elements required to generate a finitely generated subgroup is finite (locally cyclic groups are precisely the groups of rank $1$). {\it Residually finite} (\hbox{i.e.}, the finite residual is trivial) minimax groups do not have any Pr\"ufer subgroup.

\section{Dedekind skew braces and nilpotency conditions}\label{sec1vero}

\begin{defn}
{\rm Let $(B,+,\circ)$ be a skew brace. Then $B$ is said to be {\it Dedekind} if all its sub-skew braces are ideals.}
\end{defn}

\begin{rem}
{\rm The consideration of \texttt{SmallSkewbrace(8,39)} from the \texttt{YangBaxter} package\linebreak of~GAP~\cite{GAP} shows that there exist Dedekind skew braces of non-abelian type which are not {\it bi-skew braces} (that is, we cannot swap the two operations and obtain again a skew brace). Thus, one cannot hope to recover all Dedekind skew braces from the Dedekind braces discussed in~\cite{ballesterDedekindii} and \cite{dedekindbraces}. We note in passing that if $(B,+,\circ)$ is any Dedekind bi-skew brace, then $(B,\circ,+)$ is a Dedekind skew brace as well because sub-skew braces and ideals are preserved from one skew brace to the other.}
\end{rem}

\begin{rem}
 {\rm The study of skew braces is strictly related to that of the Hopf–Galois structures
on a finite Galois extension $L/K$ (see the appendix of \cite{SV}). Let $H$ be a Hopf--Galois structure on $L/K$, and $B$ the skew brace associated to $H$.  If $B$ is Dedekind, then the image of the Hopf--Galois correspondence for~$H$ is contained in the set of the sub-extensions of $L/K$ that are normal over $K$, or, equivalently, all Hopf sub-algebras $J$ of $H$ are normal and such that $L^J/K$ is~Ga\-lois and~\hbox{Hopf--Ga}\-lois with Hopf algebra $H/J$ (see \cite[Theorem~3.9, Corollary~4.1]{stefanellotrap} and~\cite[Proposition~3.4.10]{stefanello-tesi})}.   
\end{rem}

\begin{rem}
{\rm In an almost polycyclic skew brace $(B,+,\circ)$, the property of being Dedekind can be deduced from the finite homomorphic images. In fact, assume that the finite homomorphic images of $B$ are Dedekind, let $C$ be a sub-skew brace of $B$, and let $I$ be any ideal of $B$ such that~$B/I$ is finite. By hypothesis,~$B/I$ is Dedekind, so $(C+I)/I$ is an ideal of~$B/I$. Now,~The\-o\-rem~3.36 of \cite{poly} shows that $C$ is an ideal of $B$. The arbitrariness of $C$ completes the proof.}
\end{rem}


The term ‘‘Dedekind’’ comes from group theory, where a {\it Dedekind} group is one in which every subgroup is normal; non-abelian examples are called {\it Hamiltonian}. The structure of these groups was first outlined by Richard Dedekind and then completely determined by~Rein\-hold~Baer (see \cite{acourse}, 5.3.7). It turns out that a non-abelian Dedekind group is a direct product of three pieces: a quaternion group of order $8$ (which accounts for the term ‘‘Hamiltonian’’), an elementary abelian $2$-group, and a periodic abelian $2'$-group. In particular, these Dedekind groups are always nilpotent of class at most~$2$. When one begins analyzing Dedekind skew braces of small order, it quickly becomes apparent that no similarly elegant classification seems possible in the setting of skew braces (see Examples~\ref{example28},~\ref{ex29} and \ref{example41}). Indeed, in light of~\cite{ballesterDedekindii} and \cite{dedekindbraces}, even establishing a general notion of nilpotency for arbitrary skew braces appears to be a challenging task. This leads to the first problem we aim to address. Our first main result shows that, under very broad conditions (including the finite case), Dedekind skew braces are always centrally nilpotent.

\begin{theo}\label{ThA}
Let $(B,+,\circ)$ be a Dedekind skew brace having an ascending series of ideals whose factors are either locally finite or have a cyclic additive/multiplicative group.  
\begin{itemize}
    \item[\textnormal{(1)}] If $I$ is an ideal of $B$, then $I\leq Z(B)$ provided that either $(I,+)$ or $(I,\circ)$ is infinite cyclic.
    \item[\textnormal{(2)}] $B$ is hypercentral.
\end{itemize}

\noindent In particular, every trivial sub-skew brace $X$ of $B$ is contained in $Z(B)$ provided that $(X,+)$ is not periodic.
\end{theo}
\begin{proof} 
First, suppose that $B$ is finite and fix a prime $p$. We claim that~$B$ contains a sub-skew brace of order $p^k$, where $p^k$ is the maximum power of $p$ dividing $|B|$. To this end, assume that~$B$ is a minimal counterexample to the claim. By Theorem A of \cite{ballester2024soluble}, $B$ has a sub-skew brace~$I$ of prime order~$q$. By hypothesis,~$I$ is an ideal of $B$, and the minimality of $B$ implies that~$B/I$ contains a sub-skew brace $P/I$ of order the maximum power of $p$ dividing $|B/I|$. Clearly, this implies that $|P/I|=p^k$ or~\hbox{$|P/I|=p^{k-1}$.} If~\hbox{$p=q$}, then $P$ is the required sub-skew brace. Hence, we may assume that $B$ has no sub-skew brace of order $p$, and in particular that~\hbox{$p\neq q$.} Furthermore, by the minimality of $B$, we may further assume that $B/I=P/I$. According to~The\-o\-rem~A of~\cite{ballester2024soluble}, $B/I$ contains a sub-skew brace $D/I$ of order $p$. It follows that~\hbox{$|D|=pq$,} and that $D$ is supersoluble (see \cite{supersol}, Co\-rol\-la\-ry~3.9). If~\hbox{$p>q$,} then $D$ (and thus~$B$) has a sub-skew brace of order $p$ (see \cite{supersol}, The\-o\-rem~3.11), contradicting our assumptions. If~\hbox{$p<q$,} then $D$ is one of the skew braces classified in \cite{acri2020skew}, all of which have a sub-skew brace of order $p$, again contradicting our assumptions. This completes the proof of the claim.

Now, for each prime $p$, let $B_p$ be the sub-skew brace of order the maximum power of $p$ dividing the order of $B$. Since the $B_p$'s are ideals, we have that $B$ is the direct product of the~$B_p$'s. Let $I$ be a sub-skew brace of order $p$ of $B_p$.  Since $I$ is an ideal of $B$, so $I$ is contained in the centers of the~$p$-groups $(B_p,+)$, $(B_p,\circ)$ and $(B_p,+)\rtimes_\lambda(B_p,\circ)$. Consequently, $I\leq Z(B_p)$. An easy induction on the order of $B_p$ with the aid of Theorem A of \cite{ballester2024soluble} shows that $B_p$ is centrally nilpotent. It follows that $B$ is centrally nilpotent, which proves the statement when $B$ is finite.

\medskip

Assume $B$ is infinite. The condition on $B$ can be rephrased by saying that $B$ has an ascending series of ideals whose factors are cyclic groups. Since the hypothesis is inherited by homomorphic images, we only need to prove (1) for a trivial ideal $I$ (by Theorem 4.43 of~\cite{supersol}). In order to do so, we may assume that $B$ is finitely generated as a skew brace. By~\cite{trombetti2023structure},~Lem\-ma~3.1, this means we may assume that $B$ has a finite series of ideals of the previous type. Let $p$ be any prime. By induction on the Hirsch length of $(B,+)$, we have that~$B/I^p$ is centrally nilpotent, so~$I/I^p\leq Z(B/I^p)$. The arbitrariness of $p$ shows that $I\leq Z(B)$ and completes the proof of both~(1) and (2).

\medskip

Finally, if $X$ is a trivial sub-skew brace of $B$ such that $(X,+)$ is not periodic, then $X$ is a~De\-de\-kind group, so it is generated by its aperiodic elements. However, by (1), the aperiodic elements of $X$ are contained in $Z(B)$, and so $X\leq Z(B)$.
\end{proof}


\begin{rem}
{\rm Compare Theorem \ref{ThA} (1) with Example \ref{example41}.}
\end{rem}

\begin{cor}\label{cordedcentralnilp}\label{corThA}
Let $(B,+,\circ)$ be a Dedekind skew brace. If $B$ satisfies one of the following conditions $$\textit{$\bullet$\; $B$ is finite,\quad$\bullet$\; $B$ is almost polycyclic,\quad$\bullet$\; $B$ is supersoluble,}$$ then $B$ is centrally nilpotent.

\medskip

\noindent If $B$ satisfies one of the following conditions
$$\textit{$\bullet$\; $B$ is locally finite,\quad$\bullet$\; $B$ is soluble,\quad$\bullet$\; $B$ is hypercyclic,\quad$\bullet$\; $B$ is hypermultipermutational},$$ then $B$ is  hypercentral. 
\end{cor}

\begin{rem}
{\rm Corollary \ref{corThA} cannot be improved by showing that a locally finite Dedekind skew brace is centrally nilpotent. In fact, Example \ref{example28} below shows that for every odd prime $p$ there is a~De\-de\-kind skew brace $B_p$ of order $p^p$ with central nilpotency class $p$. Thus, the restricted direct product of the $B_p$'s is locally finite and Dedekind, but not centrally nilpotent (actually, not even soluble or left/right nilpotent).}
\end{rem}

\begin{cor}\label{corleftrigh}
Let $(B,+,\circ)$ be a Dedekind skew brace having a finite series of ideals whose factors are groups. \textnormal(The latter condition holds in particular when $B$ is either soluble or left/right nilpotent.\textnormal)

\begin{itemize}
    \item[\textnormal{(1)}] $B$ is soluble.
    \item[\textnormal{(2)}]  If either $\pi(B,+)$ or $\pi(B,\circ)$ is finite, and $(B,+)$ has finite Hirsch length, then $B$ is centrally nilpotent. In particular, if either $(B,+)$ or $(B,\circ)$ are $p$-groups for some prime $p$, then $B$ is centrally nilpotent.
\end{itemize}
\end{cor}
\begin{proof}
Let $\{0\}=I_0\leq I_1\leq\ldots\leq I_n=B$ be a finite series of ideals of $B$ whose factors are trivial. Note that by Theorem \ref{ThA}, $B$ is hypercentral, and so locally centrally-nilpotent. 

\medskip

\noindent(1)\quad This follows at once from the hypothesis of being Dedekind, and the characterization of~Dedekind groups.

\medskip

\noindent(2)\quad  The non-periodic factors of the series are central by Theorem \ref{ThA}, so an easy application of~Schur's theorem (see for example \cite{Ro72}, Corollary to Theorem 4.12) makes it possible to assume that, for every $i$, $(I_{i+1}/I_i,+)$ is either locally cyclic and non-periodic, or a primary group with $\pi(I_{i+1}/I_i,+)\subseteq \pi=\pi(B,+)\cup\pi(B,\circ)$. By induction, $B/I_1$ is centrally nilpotent, so we only need to show that $I_1\leq Z_\ell(B)$ for some positive integer $\ell$. Clearly, the claim holds if~$(I_1,+)$ is non-periodic. Moreover, if~$(I_1,+)$ has finite exponent, then the claim is true because $(I_1,+)$ is a Dedekind group and the subgroups of prime order are contained in~$Z(B)$ by hypercentrality of $B$. 

Assume $(I_1,+)$ is periodic of infinite exponent. Let $a$ be any element of $B$ that does not act trivially on $I_1$ with respect to either $+$, $\circ$ or $\lambda$. Since $B$ is locally centrally-nilpotent, so $\langle a\rangle$ has a finite series of ideals whose factors are cyclic groups. Consequently, $\langle a\rangle$ is contained in some finite term of the upper central series of $B$. 

First, suppose that no such an element $a$ is additively aperiodic. In this case, the three natural automorphisms groups of $(B,+)$ and $(B,\circ)$ induced on $I_1$ are finite, and the sub-skew brace $X$ generated by set of elements realizing these automorphisms is finite. Thus, $X$ is contained in a finite term of the upper central series of $B$ and $B/X$ satisfies the hypotheses. But~$I_1\circ X/X\leq Z(B/X)$, so $B/X$ is centrally nilpotent, and hence $B$ is centrally nilpotent.

Now, suppose $a$ is additively aperiodic, and it is chosen in such a way that it belongs to the smallest possible term $L$ of the series $\{I_i\}_{1\leq i\leq n}$. The previous paragraph shows that $L/\langle a\rangle$ is centrally nilpotent. Let $T/\langle a\rangle$ be the ideal of $L/\langle a\rangle$ generated by the additive/multiplicative $\pi'$-elements. Since $T$ is centrally nilpotent, we have that $T/S\leq Z(B/S)$, where $S$ is the torsion subgroup of $(T,+)$. Clearly, $S$ coincides with the torsion subgroup of $\big(\langle a\rangle,+\big)$, and so $T$ is contained in some finite term of the upper central series of $B$. Now, $B/T$ satisfies the hypotheses of (2), but has a smaller Hirsch length. Induction applies. The proof is complete.
\end{proof}

\begin{cor}\label{torsionfreeincredible}
Let $(B,+,\circ)$ be a Dedekind skew brace having an ascending series of ideals whose factors are groups \textnormal(this holds in particular when $B$ is hypermultipermutational\textnormal). If $(Z(B),+)$ is torsion-free, then $B$ is a trivial brace.
\end{cor}
\begin{proof}
By Theorem \ref{ThA}, $B$ is hypercentral. Now, Theorem 4.17 of \cite{centnilpo} shows that $\big(Z_2(B)/Z(B),+\big)$ is torsion-free. Let $b\in Z_2(B)\setminus Z(B)$. Then $\langle b, Z(B)\rangle$ is a hypermultipermutational Dedekind brace with a torsion-free additive group, so Theorem B of \cite{ballesterDedekindii} shows that $H$ is trivial. Then $\langle b\rangle=\langle b\rangle_+$ is an ideal of $B$ having trivial intersection with $Z(B)$. On the other hand, $b\in Z_2(B)$, so $[b,x]_+,[b,x]_\circ,\lambda_x(b)\in \langle b\rangle_+\cap Z(B)=\{0\}$ for all $x\in B$. Thus, $b\in Z(B)$, a contradiction. This shows that $B=Z(B)$ and proves the statement.
\end{proof}

\begin{cor}\label{solution1}
Let $(X,r)$ be a hypermultipermutational solution. If the structure skew brace $G(X,r)$ is~De\-de\-kind with a torsion-free additive group, then $(X,r)$ is a twist solution. In particular, if $(X,r)$ is a finite hypermultipermutational solution, then $(X,r)$ is a twist solution.
\end{cor}

\medskip

In the non-periodic case, the situation appears to be more difficult and the best we can do is to describe Dedekind skew braces with multipermutational level $2$.

\begin{cor}\label{classe2}
Let $(B,+,\circ)$ be a Dedekind skew brace such that $\operatorname{Soc}_2(B)=B$ and $(B,+)$ is not periodic. Then:
\begin{itemize}
    \item[\textnormal{(1)}] $\operatorname{Soc}(B)=Z(B)$, so $B=Z_2(B)$.
    \item[\textnormal{(2)}] $(B\ast B,+)$ is locally cyclic and periodic.
    \item[\textnormal{(3)}] $B$ embeds into the direct product $S\times T$ of the skew braces $S$ and $T$, where $(T,+)$ is periodic, $(T\ast T,+)$ is locally cyclic, and $S$ is a trivial skew brace with torsion-free additive group.

\end{itemize}
\end{cor}
\begin{proof}
It follows from Theorem \ref{ThA} that $B$ is locally centrally-nilpotent. Thus, Theorem 4.12 of~\cite{centnilpo} shows that the set $T$ of all periodic elements of $(B,+)$ is an ideal of $B$.

Moreover, let $M$ be a maximal ideal of $B$ such that $(M,+)$ is torsion-free. Suppose that $S/M=\operatorname{Soc}(B/M)$ is periodic, while $B/M$ is not. Let $b+M$ be an aperiodic element of $(B/M,+)$ and note that $(S+\langle b\rangle_+)/M$ is a brace, so Theorem C of \cite{ballesterDedekindii} shows that $(S+\langle b\rangle_+)/M$ has a sub-skew brace (and so an ideal of $B$) whose additive group is torsion-free. The maximality of $M$ gives a contradiction. Thus, $(B/M,+)$ is periodic. 

Therefore, $B$ embeds into $B/M\times B/T$. It follows from Corollary \ref{torsionfreeincredible} that $B/T$ is trivial. Therefore the socle of $B/M\times B/T$ contains $B/T$, and this implies that $(B/\operatorname{Soc}(B),+)$ is periodic, and that $(\operatorname{Soc}(B),+)$ is not. In particular, Theorem \ref{ThA} shows that  $\operatorname{Soc}(B)=Z(B)$.

Now, let $b\in B$ be of infinite additive order. By what we have already proved, $b\ast b$ is additively periodic, so $\langle b\ast b\rangle_+$ is a finite ideal of $B$. Moreover, the sub-skew brace $U/\langle b\ast b\rangle_+$ generated by $b+\langle b\ast b\rangle_+$ is trivial by Lemma \ref{bastb}, and so is an infinite cyclic group. Since~$B$ is locally centrally-nilpotent, so $U/\langle b\ast b\rangle_+$ is contained in $Z\big(U/\langle b\ast b\rangle_+\big)$. Consequently, $x\ast b$ and $b\ast x$ belong to $\langle b\ast b\rangle_+$ for every $x\in B$.

Let $b\in B$ be of finite additive order, and let $d$ be an additively aperiodic element of $Z(B)$. Then $x=b+d$ has infinite additive order, so $b\ast x$ and $x\ast b$ belong to $\langle (d+b)\ast (d+b)\rangle_+=\langle b\ast b\rangle_+$ for every $x\in B$. Therefore in any case $b\ast x$, and $x\ast b$ belong to $\langle b\ast b\rangle_+$ for every $b,x\in B$.

Now, we show that $B\ast B$ is locally cyclic. To this aim, let  $a,b\in B$. We claim that either $\langle a\ast a\rangle_+\subseteq\langle b\ast b\rangle_+$ or $\langle b\ast b\rangle_+\subseteq\langle a\ast a\rangle_+$. Since $B$ is locally centrally-nilpotent, we may assume without loss of generality that $\langle a,b\rangle$ is a finite sub-skew brace of order a power of a prime $p$. Of course, we may further assume $a\ast a\neq0\neq b\ast b$.

Suppose first $\langle a\ast a\rangle_+\cap\langle b\ast b\rangle_+=\{0\}$. In particular, $\langle a\ast a,b\ast b\rangle_+=\langle a\ast a\rangle_+\times\langle b\ast b\rangle_+$ is the direct product of two non-zero cyclic groups of order a power of $p$. Since $a\ast b,b\ast a\in \langle a\ast a\rangle_+\cap\langle b\ast b\rangle_+=\{0\}$, so $$
\begin{array}{c}
(a+b)\ast (a+b)=(a\circ b)\ast (a+b)=(a\circ b)\ast a+(a\circ b)\ast b\\[0.2cm]
=a\ast (b\ast a)+a\ast a+b\ast a+a\ast (b\ast b)+a\ast b+b\ast b=a\ast a+b\ast b,
\end{array}
$$ and hence $$a\ast a=a\ast (a+b)\in\langle (a+b)\ast (a+b)\rangle_+=\langle a\ast a+b\ast b\rangle_+.$$   However, this is clearly impossible because $\langle a\ast a+b\ast b\rangle_+$ would have to be not cyclic in this case. Thus, $V=\langle a\ast a\rangle_+\cap\langle b\ast b\rangle_+\neq\{0\}$. Now, the previous argument applies to $B/V$, and hence either $a\ast a$ or $b\ast b$ must be contained in $V$, which means that either $\langle a\ast a\rangle_+\leq\langle b\ast b\rangle_+$ or $\langle b\ast b\rangle_+\leq\langle a\ast a\rangle_+$. This implies that $(B\ast B,+)$ is locally cyclic.    
\end{proof}

\medskip

We finish this section with some natural open questions and relevant examples.

\begin{quest}
Let $(B,+,\circ)$ be a Dedekind skew brace.
\begin{itemize}
    \item Is $B$ hypercentral? 
    \item Does the hypersocle coincide with the hypercentre? \textnormal(at least under the hypothesis that the additive group is torsion-free\textnormal)
    \item Can the hypothesis on the Hirsch length  in \textnormal{Corollary \ref{corleftrigh} (2)} be dropped?
\end{itemize}
\end{quest}

\begin{ex}\label{example28}
For each prime $p$ and any positive integer $n\geq2$, there exists a Dedekind bi-skew brace of central nilpotency class $n$.
\end{ex}
\begin{proof}
Suppose first $p=2$. If $n=2$, then \texttt{SmallSkewbrace(4,2)} is the required skew brace. Suppose $n\geq3$, and put $(B,+)=(\mathbb{Z}_{2^n},+)$. Define $a\circ b=a+(-1)^ab$, for all $a,b\in B$. It is easy to see that $(B,+,\circ)$ is a bi-skew brace (see \cite{caranti}, Theorem 3.1 (3)). The fact that $B$ is~De\-de\-kind can either be checked by direct computation or by using Corollary \ref{corfinitecycliccase}. Of course, the nilpotency class of $(B,\circ)$ is always $n-1$, so also the central nilpotency class of $B$ is $n-1$ (note however that $\operatorname{Soc}_2(B)=B$).

Assume $p>2$ and let $(B,+)=\big(\langle x\rangle,+\big)$ be a cyclic group of order $p^n$. 
By \cite[Proposition~2.19]{campedel} the assignment $\lambda_x(x)=(p+1)x$ defines a unique map
$$\lambda\colon B\mapsto\operatorname{Aut}(B,+)$$ verifying the equality 
$$\lambda_{y+\lambda_y(z)}=\lambda_y\lambda_z\qquad\forall y,z\in B$$
(in the cited paper called ``gamma function equation'')
and therefore it also defines an operation $\circ$ on $B$ such that 
$$x^\ell=\underbrace{x\circ x\circ\ldots\circ x}_{\ell\textnormal{ times}}=
\big(1+(1+p)+\ldots+(1+p)^{\ell-1}\big)x.
$$
(see \cite{campedel} after Theorem 2.2, or \cite[Section 2]{idc2023} for definitions given with the left notation). 
By~Corollary~2.18 of \cite{campedel}, $(B,\circ)$ is cyclic and generated by $x$. The fact that $(B,+,\circ)$ is a bi-skew brace follows at once from Lemma 1.10 of \cite{pompili} and~The\-o\-rem~3.1 (3) of \cite{caranti}, while the fact that it is Dedekind can either be checked by direct computation or by using Corollary \ref{corfinitecycliccase}.




Now, the centre of $B$ is $K=\operatorname{Ker}(\lambda)$, which has order $p$. Since $(B/K,+,\circ)$ is isomorphic to the example we just constructed but with order $p^{n-1}$, so $B$ has central nilpotency class $n$ (actually,~$B$ has multipermutational level $n$).
\end{proof}

\begin{ex}\label{ex29}
For each odd prime $p$ and any positive integer $n\geq2$ there exists a \textnormal(non-trivial\textnormal)~De\-de\-kind bi-skew brace $B$ of order $p^{n+1}$ whose additive and multiplicative groups are non-abelian $p$-groups of class~$2$. Moreover, $B$ is centrally nilpotent of class $2$.
\end{ex}
\begin{proof}
Choose $n\geq2$ and define $(B,+)=\langle y\rangle\ltimes\langle x\rangle$, where $o(x)=p^n$, $o(y)=p$ and $-y+x+y=(1+p^{n-1})x$. Then $(B,+)$ has order $p^{n+1}$ and class $2$. Now, the conjugacy by $x$ defines an automorphism of order $p$ because $$[px,y]=p[x,y]=p^nx=0,$$ and so the assignation $$
\begin{cases}
x\mapsto \operatorname{Id}_{(B,+)}\\
y\mapsto \big(\iota(x): \ast\mapsto -x+\ast+x\big)
\end{cases}
$$ defines a morphism $\lambda: b\in (B,+)\mapsto \lambda_b\in\operatorname{Aut}(B,+)$. In particular, $\lambda_{\lambda_b(c)}=\lambda_{c}$ for eve\-ry~\hbox{$c\in (B,+)$,} and so the position $a\circ b=a+\lambda_a(b)$ defines a skew brace by Lemma 2.13 of~\cite{campedel}; this is a bi-skew brace by \cite{caranti}, Theorem 3.1 (3).

Note that $$
\begin{array}{c}
x^{-1}\circ y^{-1}\circ x\circ y
=-x+y^{-1}\circ x\circ y
=-x+y^{-1}\circ (x+y)
=-x+y^{-1}\circ x-y^{-1}\\[0.2cm]
=-x+y^{-1}+\lambda_{y}^{-1}(x)-y^{-1}
=-x+y^{-1}+x-y^{-1}
=-x+ x-y-x +x+x+y-x\\[0.2cm]
=-y+x+y-x=-x-y+x+y,
\end{array}
$$ and $$
\begin{array}{c}
y^2=y\circ y=y-x+y+x=2y-p^{n-1}x,\\[0.2cm]
\ldots\ldots\ldots\ldots,\\[0.2cm]
y^k=y\circ y^{k-1}=y-x+y^{k-1}+x
=y-x+(k-1)y-\binom{k-1}{2}p^{n-1}x+x\\[0.2cm]
=ky-(k-1)p^{n-1}x-\binom{k-1}{2}p^{n-1}x=ky-\binom{k}{2}p^{n-1}x.
\end{array}
$$
The previous computations show that $(B,\circ)\simeq (B,+)$. 

Now we prove that $(B,+,\circ)$ is Dedekind. Let $H$ be any sub-skew brace of $B$. If $H$ is contained in $\langle x\rangle_+$, then $(H,+)$ and $(H,\circ)$ are characteristic in $\big(\langle x\rangle_+,+\big)$ and $\big(\langle x\rangle_+,\circ\big)$, respectively, so $H$ is an ideal of $B$. Suppose $H$ is not contained in $\langle x\rangle_+$. Then the additive subgroup of $H$ must contain $U=\langle p^{n-1}x\rangle_+$. It is pretty clear that $U=Z(B)$, so $H/U$ is an ideal of $B/U$ and so $H$ is an ideal of $B$. The proof is complete.
\end{proof}

\begin{ex}\label{example41}
Let $p$ be a prime. There exist non-trivial infinite supersoluble Dedekind skew braces containing elements of additive order $p$.
\end{ex}
\begin{proof}
Let $p$ be an odd prime, and let $A$ be any non-trivial skew brace of order $p^2$ with an elementary abelian additive group. It can be easily verified from the description in \cite{p3} that~$A$ has only one sub-skew brace of order $p$, which is in fact an ideal. Thus,~$A$ is Dedekind. Let $B$ be the direct
product of $A$ and an infinite cyclic group~$C$, and let $b\in B$ such that $b=a+c=a\circ c$ for some~\hbox{$a\in A\setminus\{0\}$} and $c\in C\setminus \{0\}$. Put~\hbox{$U=\langle b\rangle$.} Then~$\langle c^p\rangle_\circ$ is contained in $U$, so we may assume $C$ has order~$p$. Now, if~\hbox{$a\in Z(B)$,} then $b\in Z(B)$ and~$U$ is an ideal. Assume $a\notin Z(B)$. If $\lambda_a(a)\in\langle a\rangle_+$, then~$\langle a\rangle_+$ is a sub-skew brace of order~$p$ which is different from $Z(A)$. Moreover, since $a\not\in Z(A)$, we have that~$\langle a\rangle_+$ is~\hbox{$\lambda$-inva}\-riant and so an ideal of $A$, which implies that $A=Z(A)\times \langle a\rangle_+$ is abelian, a contradiction. Thus, $\lambda_a(a)\notin \langle a\rangle_+$. Write $\lambda_a(a)=a+d$, where $d\in Z(A)\setminus\{0\}$. Consequently,~$U$ contains $\lambda_{a\circ c}(a+c)=a+d+c$ and so also $d\in Z(A)$. Finally, since $U/Z(A)$ is an ideal of the trivial brace $B/Z(A)$, so $U$ is an ideal of~$B$. 

A similar example also works for $p=2$, but in this case there are a couple of additional possibilities. Let $A$ be a skew brace of order $4$ in which the additive group and the multiplicative group are not isomorphic; for the sake of simplicity, assume $(A,+)$ is cyclic and $(A,\circ)$ is elementary abelian. It follows from Theorem \ref{thB} that $A$ is Dedekind. Let $B$ be the direct product of~$A$ and an infinite cyclic group $C$. Let $b\in B$ and write $b=a+c=a\circ c$ for some~\hbox{$a\in A\setminus\{0\}$} and $c\in C\setminus \{0\}$. We claim that the sub-skew brace $U$ generated by $b$ is an ideal of $B$. Clearly,~\hbox{$c^2\in U$,} so factoring by $\langle c^2\rangle_\circ$, we may assume $C$ is finite of order $2$. Now, if $a$ has additive order $2$, then $b\in Z(B)$, and we are done. If $a$ has additive order $4$, then $\langle a\rangle_+=A$, so $U$ contains the ideal~$(2A,+)$, and again $U$ is an ideal because $B/2A$ is a trivial brace. The arbitrariness of $b$ shows that $B$ is Dedekind.
\end{proof}

\section{Locally cyclic groups in Dedekind skew braces}

The aim of this section is to study Dedekind skew braces in certain finite rank environments. We start by investigating the case in which one of the two groups is finite cyclic. Our first result extends Proposition~11 of \cite{dedekindbraces} and gives us a sufficient condition for establishing that a skew brace of prime power order is~Dedekind.

\begin{theo}\label{thB}
Let $p$ be a prime, and $(B,+,\circ)$ a finite skew brace of order $p^n$. If either $(B,+)$ or $(B,\circ)$ is cyclic, then $B$ is a Dedekind skew brace. In particular, $B$ is centrally nilpotent.
\end{theo}
\begin{proof}
First, we observe that if $B$ has a sub-skew brace of order $p^k$, then by hypothesis $B$ has precisely one sub-skew brace of this order.

Now, since the natural semidirect product $(B,+)\rtimes_\lambda(B,\circ)$ is nilpotent, there exists a strong left-ideal $I$ of $B$ of order $p$. If $(B,\circ)$ is abelian, then $I$ is an ideal of $B$. Assume~$(B,\circ)$ is not abelian, so $(B,+)$ is cyclic. Consequently, $(B,\circ)/\operatorname{Ker}(\lambda)$ is abelian, and \hbox{$\operatorname{Ker}(\lambda)\neq\{0\}$.} Since~$\operatorname{Ker}(\lambda)$ is a proper non-zero sub-skew brace of $B$, so it is Dedekind by induction, and hence it must contain a sub-skew brace of order $p$ (see Theorem \ref{ThA}). In turn, this latter sub-skew brace must coincide with $I$ (by our starting remark), so in particular, $I\ast B=\{0\}$. It follows from~Lem\-ma~\ref{lemuseful} that $I$ is an ideal of $B$ in any case.

Finally, if $C$ is any proper sub-skew brace of $B$, then $C$ is Dedekind by induction so it contains the ideal $I$ of order $p$ by the previous paragraph. By induction, $B/I$ is Dedekind, so~$C/I$ is an ideal of $B/I$. Consequently, $C$ is an ideal of $B$, and the statement is proved.
\end{proof}

\begin{cor}\label{corfinitecycliccase}
Let $(B,+,\circ)$ be a finite skew brace such that either $(B,+)$ or $(B,\circ)$ is cyclic and the other nilpotent. Then $B$ is Dedekind.
\end{cor}
\begin{proof}
In either case, the Sylow subgroups of $(B,+)$ are ideals and we may assume that \hbox{$|B|=p^n$} for some prime $p$. Then the statement follows from Theorem \ref{thB}.
\end{proof}

\medskip

Theorem \ref{thB} cannot be much improved because there exist non-Dedekind skew braces of composite order in which either the additive or the multiplicative group is cyclic. Moreover, there exist non-Dedekind skew braces with abelian additive and multiplicative groups (for example \texttt{SmallBrace(8,7)}). Despite these examples, we can clarify a bit better what happens in case one of the two groups is finite cyclic. First, we give an easy proof (using Dedekind-ness) of a result that was essentially already known in the Hopf--Galois context (see 
\cite{campedel}, Co\-rol\-la\-ry~2.18, and \cite{Kohl}, Theorem 4.5).

\begin{theo}\label{cyclicideals}
Let $p$ be an odd prime and $(B,+,\circ)$ a skew brace of finite order $p^n$. If either $(B,+)$ or~$(B,\circ)$ is cyclic, then the other is cyclic as well. Moreover, additive and multiplicative generators coincide.
\end{theo}
\begin{proof}
By Theorem \ref{thB}, $B$ is centrally nilpotent, so it has an ideal $I$ of order $p$; in particular,~\hbox{$I\leq Z(B)$.} By induction, $(B/I,+)$ and $(B/I,\circ)$ are cyclic with the same generators. Being central-by-cyclic, both $(B,+)$ and $(B,\circ)$ are abelian. Now, if either $(B,+)$ or~$(B,\circ)$ is not cyclic, then the only ideal $X/I$ of order~$p$ of $B/I$ is such that either $(X,+)$ or $(X,\circ)$ is not cyclic. But~$X$ has order $p^2$, and so the characterization of skew braces of order $p^2$ given in \cite{p3} gives a contradiction. This also implicitly implies that the additive and multiplicative generators are the same.
\end{proof}

\begin{rem}
{\rm Theorem \ref{cyclicideals} does not hold for skew braces whose order is not a power of an odd prime. This is shown for example in \cite{p3}, but you may also easily find on GAP examples of skew braces in which either the additive or the multiplicative group is cyclic of  composite order but the other is not (see also \cite{byott}).}
\end{rem}

\begin{theo}\label{cyclicsupersoluble}
Let $(B,+,\circ)$ be a finite skew brace. If either $(B,+)$ or $(B,\circ)$ are cyclic, then $B$ is supersoluble.
\end{theo}
\begin{proof}
If the order of $B$ is a power of a prime, then the result follows from~The\-o\-rem~\ref{cyclicideals}. Thus, we may assume the order of $B$ is not a prime power.



\medskip

First, suppose that $(B,+)$ is cyclic, and let $p$ be the smallest prime dividing the order of~$B$. Since $(B,+)$ is cyclic, there is a characteristic subgroup $(C,+)$ of index $p$ in $(B,+)$. In particular,~$C$ is a sub-skew brace of index $p$. Since $(C,\circ)$ has index $p$ in $(B,\circ)$, so $(C,\circ)$ is multiplicatively normal, and hence $C$ is an ideal of $B$. By induction on the order, $C$ is supersoluble. Let~$q$ be the largest prime dividing the order of~$B$ (note that $q>p\geq2$), and let $(E,+)$ be a~Sy\-low~\hbox{$q$-sub}\-group of $(B,+)$. Then $E$ is an ideal of $B$. Let $Q$ be the only additive and multiplicative subgroup of order $q$ of $E$ (see Theorem \ref{cyclicsupersoluble}). Then $Q$ is an ideal of~$B$, and induction yields that~$B/Q$ is supersoluble, so $B$ is supersoluble.

\medskip

Next, suppose $(B,\circ)$ is cyclic. Since $B$ is two-sided, $(B\ast B,+)$ is abelian (see Lemma 4.5 of~\cite{nasy}). The result is obvious if $B\ast B=\{0\}$, so we may assume this is not the case. Then the~Sy\-low~sub\-groups of $(B\ast B,+)$ are ideals of $B$. Let $P$ be one of these ideals of order $p^n$. If $p>2$, then~$(P,+)$ is cyclic (by Theorem \ref{cyclicideals}) and consequently the only additive and multiplicative subgroup of order~$p$ of $P$ is an ideal of $B$. Induction completes the proof in this case. 

Assume $B\ast B$ has order a power of $2$ (so it is centrally nilpotent by Theorem \ref{thB}). We want to find an ideal of odd prime order of $B$, so we can then apply induction. In order to do so, we may assume $B/B\ast B$ has odd prime order (recall that $B/B\ast B$ is supersoluble) and that $(B,+)$ is not abelian. Since $\big(2(B\ast B),+\big)$ is an ideal we may further assume $(B\ast B,+)$ is elementary abelian. Now, one can easily check on GAP (see also \cite{p3}) that there are no skew braces of order~$8$ whose additive group is elementary abelian and the multiplicative is cyclic. Consequently, $B\ast B$ must have order at most $4$, so $B/B\ast B$ must have order $3$. A further check on GAP of the skew braces of order $6$ and $12$ reveals that they all have an ideal of order $3$, thus giving the contradiction that $(B\ast B,+)$ is abelian.
\end{proof}


\begin{rem}
{\rm The condition in Theorem \ref{cyclicsupersoluble} that one of the groups is cyclic is mandatory because there exist finite skew braces of prime power order that are not centrally nilpotent (so they cannot be supersoluble by \cite{supersol}, Theorem 3.7).}
\end{rem}

\medskip

Next, we move to the infinite case.



%

\begin{theo}\label{charQ+}
Let $(B,+,\circ)$ be a non-trivial skew brace in which $(B,+)$ is locally cyclic and torsion-free. Let $K=\operatorname{Ker}(\lambda)$ and consider $(B,+)$ as a subgroup of $(\mathbb Q,+)$ containing the rational number $1$. Then one of the following possibilities holds\textnormal:
    \begin{itemize}
        \item[\textnormal{(1a)}] $K$ is a non-zero ideal of $B$ such that $|B/K|=2$, $(B,+)$ is not $2$-divisible, and 
        $$a\circ b=a+(-1)^{\varphi(a+K)}b$$ for all $a,b\in B$, where $\varphi:(B/K,+)\rightarrow\mathbb{Z}_{2}$ is an isomorphism. Thus, $(B,\circ)$ is isomorphic to the dihedral group $\mathbb Z_2\ltimes (B,+)$. In this case, $B$ is not Dedekind.


\item[\textnormal{(1b)}] $K=\{0\}$, $(B,\circ)$ is abelian and there exists a non-zero rational number $\frac{m_1}{m_2}\in B$ such that:
\begin{itemize}
    \item[\textnormal{(i)}] $(m_1,m_2)=1$, $m_2>0$, and $m_2-m_1\neq0,\pm1$;
    \item[\textnormal{(ii)}] $(B,+)$ is $p$-divisible for every prime $p$ that does not divide $m_1-m_2$;
    \item[\textnormal{(iii)}] $(B,+)$ is not $(m_2-m_1)$-divisible;
    \item[\textnormal{(iv)}] $\frac{-a}{1-a+\frac{m_{1}}{m_{2}}a}\in B$, $(m_2-m_1)ab\in B$, and $$a\circ b=a+b-ab+\frac{m_{1}}{m_{2}}ab$$ for all $a,b\in B$.
\end{itemize}  
 In this case, $(B,+)$ is actually isomorphic to a sub-ring of $(\mathbb Q,+)$. Moreover, $B$ is not Dedekind.

    \end{itemize} 

\noindent Conversely, let $(X,+)$ be any subgroup of $(\mathbb Q,+)$ containing the rational number $1$.  
\begin{itemize}
    \item[\textnormal{(2a)}] If $X$ is not~\hbox{$2$-di}\-vi\-sible, then we can define a group $(X,\circ_1)$ by putting $$x\circ_1 y=x+(-1)^{\varphi(x+2X)}y$$ for every $x,y\in X$, where $\varphi:X/2X\rightarrow \mathbb Z_2$ is an isomorphism. Then $(X,+,\circ_1)$ defines a skew brace such that $\operatorname{Ker}(\lambda)=2X\neq\{0\}$. 
    

    \item[\textnormal{(2b)}] If $X$ is a sub-ring of $\mathbb Q$ for which there exists a non-zero rational number $\frac{m_1}{m_2}$ satisfying: 
    \begin{itemize}
    \item[\textnormal{(i)}] $(m_1,m_2)=1$, $m_2>0$, and $m_2-m_1\neq0,\pm1$;
    \item[\textnormal{(ii)}] $(X,+)$ is $p$-divisible for every prime $p$ that does not divide $m_1-m_2$;
    \item[\textnormal{(iii)}] $(B,+)$ is not $(m_2-m_1)$-divisible.
\end{itemize} 
\noindent Then the operation $$x\circ_2y=x+y-xy+\frac{m_1}{m_2}xy\qquad (x,y\in X)$$ defines an abelian group $(X,\circ_2)$, and $(X,+,\circ_2)$ defines a skew brace with $\operatorname{Ker}(\lambda)=\{0\}$.
\end{itemize}
\end{theo}
\begin{proof}
It follows from Lemma \ref{lemuseful} that $K$ is an ideal of $B$. Suppose first that $K\neq\{0\}$. Then $(B/K,+)$ is periodic and locally cyclic, so $B/K$ is covered by finite left ideals, and hence $(B/K,\circ)$ is periodic too. It follows that $(B/K,\circ)=(B/K,+)$ has order $2$. Let $\varphi:(B/K,+)\rightarrow\mathbb{Z}_{2}$ be an isomorphism. Then $$a\circ b=a+(-1)^{\varphi(a+K)}b$$ for all $a,b\in B$. 

Take $x\in B\setminus K$ and $y\in K\setminus\{0\}$, and suppose $(B,\circ)$ is abelian. Then
$$x-y=x\circ y=y\circ x=y+x=x+y\implies y=0,$$ an obvious contradiction. Consequently, $(B,\circ)$ is non-abelian and $x$ must act as the inversion on the torsion-free locally cyclic group $(K,\circ)=(K,+)\simeq (B,+)$. In order to complete the proof of (1a), we just need to observe that $(3K,+)=(3K,\circ)$ is an ideal of $B$ such that~$(B/3K,\circ)$ is not nilpotent, so Theorem \ref{ThA} yields that $B/3K$ is not Dedekind. Consequently,~$B$ is not Dedekind and (1a) is proved.


Assume conversely that $(X,+)$ is a subgroup of $(\mathbb Q,+)$ that is not $2$-divisible. For each $x\in X$, let   
$\lambda\colon X\to \operatorname{Aut}(X,+) $  defined by  $x\mapsto \lambda_x=(-1)^{\varphi(x+2X)}\operatorname{id}$. Clearly, $\operatorname{Ker}(\lambda)=2X$. To show that~$\lambda$ defines a skew brace structure on $X$
 it is enough to notice that, for $x,y\in X$,
 $$\lambda_{x+\lambda_x(y)}=\lambda_x\lambda_y$$
 and this is true since both members are the equal to $-\operatorname{id}$ when exactly one between $x$ and $y$ belongs to $2X$ and is the identity otherwise.

\medskip

Now, suppose that $K=\{0\}$, so in particular $(B,\circ)$ is abelian, being isomorphic to a subgroup of $\operatorname{Aut}(\mathbb Q,+)$. Since $\lambda_1$ is an automorphism of $(B,+)$, so there are coprime integers~$m_{1}\neq0$ and \hbox{$m_{2}>0$} such that $\lambda_{1}(x)=\frac{m_{1}}{m_{2}}x$ for every $x\in B$. Similarly, for every $a\in B$, there are coprime integers~$n_{1,a}\neq0$ and $n_{2,a}>0$ such that $\lambda_{a}(x)=\frac{n_{1,a}}{n_{2,a}}x$ for all $x\in B$. We have that
    $$a+\frac{n_{1,a}}{n_{2,a}}=a+\lambda_{a}(1)=a\circ1=1\circ a=1+\lambda_{1}(a)=1+\frac{m_{1}}{m_{2}}a,$$ which implies that $$\frac{n_{1,a}}{n_{2,a}}=1-a+\frac{m_1}{m_{2}}a$$ for every $a\in B$. Therefore $$a\circ b=a+\lambda_{a}(b)=a+b-ab+\frac{m_{1}}{m_{2}}ab$$ for every $a,b\in B$. In particular, $\big(\frac{m_1}{m_2}-1\big)ab\in B$ for every $a,b\in B$.

    Now, the fact that $n_{1,a}\neq0$ implies that the rational number $\frac{1}{1-\frac{m_1}{m_2}}=\frac{m_2}{m_2-m_1}$ does not belong to $B$, so in particular $m_2-m_1$ does not divide $m_2$ and hence $m_2-m_1\neq\pm1$. Note also that if $m_2=m_1$, then $a\circ b=a+b$ for all $a,b\in B$, so $B$ is trivial, contradicting the hypothesis.

We show that $(B,+)$ is $p$-divisible for every prime $p$ that does not divide $m_1-m_2$. First, note that if $n$ is any integer, then the $\circ$-inverse of $nm_2$ is $\frac{-nm_2}{1+(m_1-m_2)n}$. Since $p$ does not divide $m_1-m_2$, Bézout Lemma shows that for every choice of the positive integer $u$, there are integers $v,w$ such that $1+(m_1-m_2)v=p^uw$. Since  $p^uw$ is coprime with $v$, so 
$$\frac{\langle m_2\rangle_+}{\langle vm_2\rangle_+}\cap\frac{\langle vm_2/p^uw\rangle_+}{\langle vm_2\rangle_+}=\big\{\langle vm_2\rangle_+\big\}.$$ Consequently, $(B/\langle m_2\rangle_+,+)$ contains a Pr\"ufer $p$-subgroup, and so also $(B/\langle 1\rangle,+)$ has such a subgroup. This implies that $(B,+)$ is $p$-divisible.

Since $m_2$ does not divide $m_1-m_2$, the condition $\frac{m_2}{m_2-m_1}\not\in B$ can be rephrased as the condition \hbox{$\frac{1}{m_2-m_1}\not\in B$} by using $m_2$-divisibility. Similarly, the condition $\big(\frac{m_1}{m_2}-1\big)ab\in B$ for every $a,b\in B$, can be rephrased as the condition $(m_2-m_1)ab\in B$. Finally, the $\circ$-inverse of~\hbox{$a\in B$} is $\frac{-a}{1-a+\frac{m_1}{m_2}a}$, which then obviously belongs to $B$. The fact that $(B,+)$ is isomorphic to a sub-ring of $(\mathbb Q,+)$ depends on the well-known characterization of the additive subgroups of~$(\mathbb Q,+)$.

In order to complete the proof of (1b), we need to prove that $B$ is not Dedekind. Let $p$ be a prime number that does not divide $m_2(m_1-m_2)$, and let $\big(X/\mathbb Z,+\big)=O_{p'}(B/\mathbb Z,+)$; clearly,~$X$ is a sub-ring of $\mathbb Q$. Put $Y=Xp=\{xp\,:\, x\in X\}$. Then $Y$ is a proper subgroup of $(B,+)$, which is not $p$-divisible. The rule for multiplication shows that $Y$ is multiplicatively closed. Let $a=\frac{np}{m}\in Y$, where $(np,m)=1$. The multiplicative inverse of $a$ is $\frac{-npm_2}{mm_2+(m_1-m_2)np}$, which still belongs to $Y$. Thus, $Y$ is a sub-skew brace. Now, $\lambda_{1/p^2}(p)=\frac{p^2m_2-m_2+m_1}{m_2p^2}\cdot p=\frac{p^2m_2-m_2+m_1}{m_2p}$ does not belong to $Y$, so $Y$ is not an ideal and $B$ is not Dedekind.

\medskip

Suppose conversely that $X$ contains a rational number $\frac{m_1}{m_2}$ satisfying the hypotheses of (2b). These hypotheses shows that the operation $$x\circ_2y=x+y-xy+\frac{m_1}{m_2}xy\qquad (x,y\in X)$$ is well-defined, commutative, and has the rational number $0$ as the identity element.

In order to show that every element of $X$ is invertible with respect to $\circ_2$, we need to show that $\frac{-x}{1-x+\frac{m_{1}}{m_{2}}x}\in X$ for every $x\in X$. First, $1-x+\frac{m_{1}}{m_{2}}x\neq0$ because $\frac{m_2}{m_2-m_1}\not\in B$ by~(ii) and~(iii). Now, write $x=\frac{n}{m}$, where $n,m$ are integers such that $(n,m)=1$ and $m>0$. By using (ii), the condition $\frac{-x}{1-x+\frac{m_{1}}{m_{2}}x}\in X$ can be rephrased as $\frac{n}{m_2m+(m_1-m_2)n}\in X$. We need to understand if the reduced form of the fraction $\frac{n}{m_2m+(m_1-m_2)n}$ belongs to $X$. In order to do this, by (ii), we may further assume that $(m_2,n)=1$, so the condition can be rephrased as $\frac{1}{m_2m+(m_1-m_2)n}\in X$. This shows that if $p$ is any prime dividing the latter denominator and $m_2-m_1$, then $p$ must also divide $m$. Thus,  every prime divisor $q$ of $m_2m+(m_1-m_2)n$ is such that $X$ is $q$-divisible, and hence the condition $\frac{1}{m_2m+(m_1-m_2)n}\in X$ is verified.

Moreover, $$
\begin{array}{c}
x\circ_2(y\circ_2z)=x\circ_2\big(y+z+\big(\frac{m_1}{m_2}-1\big)yz\big)\\[0.2cm]
=x+y+z+\big(\frac{m_1}{m_2}-1\big)yz+\big(\frac{m_1}{m_2}-1\big)xy+\big(\frac{m_1}{m_2}-1\big)xz+\big(\frac{m_1}{m_2}-1\big)^2xyz\\[0.2cm]
=
x+y+\big(\frac{m_1}{m_2}-1\big)xy+z+\big(\frac{m_1}{m_2}-1\big)xz+\big(\frac{m_1}{m_2}-1\big)yz+\big(\frac{m_1}{m_2}-1\big)^2xyz\\[0.2cm]
=
\big(x+y+\big(\frac{m_1}{m_2}-1\big)xy\big)\circ_2z=(x\circ_2y)\circ_2z,
\end{array}
$$ for all $x,y,z\in X$. Thus, $(B,\circ_2)$ defines an abelian group.

Finally, $(X,+,\circ_2)$ is a skew brace because $$
\begin{array}{c}
x\circ_2(y+z)=x+y+z+\big(\frac{m_{1}}{m_{2}}-1\big)x(y+z)
=x+y-x+x+z+\big(\frac{m_{1}}{m_{2}}-1\big)xy+\big(\frac{m_{1}}{m_{2}}-1\big)xz\\[0.2cm]
=x+y+\big(\frac{m_{1}}{m_{2}}-1\big)xy-x+x+z+\big(\frac{m_{1}}{m_{2}}-1\big)xz=x\circ_2 y-x+x\circ_2 z,
\end{array}
$$ for all $x,y,z\in X$. Therefore (2b) is proved.
\end{proof}

\begin{rem}
{\rm The skew braces $(B,+,\circ)$ described in Theorem \ref{charQ+} (1a) and (1b) are residually centrally-nil\-po\-tent. In fact, in both cases, there is a prime $p$ for which $(B,+)$ is not $p$-divisible, and $(p^nB,+)$ is an ideal of $B$ for every positive integer $n$. Then $(B/p^nB,+,\circ)$ is a finite~De\-de\-kind skew brace, and hence centrally nilpotent by Corollary \ref{corThA}. Clearly, $\bigcap_{n\in\mathbb N}p^nB=\{0\}$ and so $B$ is residually centrally-nilpotent.}
\end{rem}

\begin{cor}
    If $(B,+,\circ)$ is a skew brace such that $(B,+)\simeq (\mathbb{Q},+)$, then $B$ is trivial.
\end{cor}

\begin{lem}\label{lemOld4.10}
Let $(X, \circ)$ be an abelian group which is not $2$-divisible, $A$ a subgroup of $(X,\circ)$ of index~$2$, and $x \in X\setminus A$. Let $\lambda$ be the morphism from $(X,\circ)$ to $\operatorname{Aut}(X,\circ)$ that
  has $A$ as kernel and sends $x$ to the inversion map $u \mapsto
  u^{-1}$ on $X$.   Then\textnormal:
  \begin{enumerate}
  \item[\textnormal{(1)}]
    $\lambda$ defines a skew brace structure
    $\muladd{X}$ via $u + v = u \circ \lambda_{u}(v)$;
  \item[\textnormal{(2)}] 
    the group $\add{X}$ is the semidirect product of a cyclic
    group of order $2$ generated by $x$ by $A$, with~$x$ acting on $A$ as the inversion;
  \item[\textnormal{(3)}]
    $(X,\circ,+)$ is a bi-skew brace, so $(X,+,\circ)$ is a skew brace too.
  \end{enumerate}
\end{lem}

\begin{proof}
  For $a, b \in A$ and $i, j = 0, 1$ we have that
  \begin{align*}
    b^{-1} \circ x^{-j}
    \circ
    \lambda_{a \circ x^{i}}(b \circ x^{j})
    &=
    b^{-1} \circ x^{-j}
    \circ
    b^{(-1)^{i}} \circ x^{(-1)^{i} j}
    \\&=
    \begin{cases}
      1 & \text{if $i = 0$}\\
      b^{-2} \circ x^{- 2 j} & \text{if $i = 1$,}
    \end{cases}
  \end{align*}
  is an element of $A$. By~\cite[Lemma 2.13]{campedel}, $\lambda$ defines a skew brace structure
    $(X,\circ,+)$ via $u + v = u \circ \lambda_{u}(v)$. Moreover, we have $x + x = x \circ \lambda_{x}(x) = x \circ x^{-1} = 1$, and clearly
  $(A,+) = (A,\circ)$. For $a \in A$ we have
  \begin{equation*}
    - x + a + x
    =
    (x \circ a^{-1}) + x
    =
    x \circ a^{-1} \circ x^{-1}
    =
    a^{-1}
    =
    -a.
  \end{equation*}

  Since $\lambda : \mul{X} \to \operatorname{Aut}(X,\circ)$ is a morphism with an
  abelian (cyclic) image, $(X,\circ,+)$ is a bi-skew brace by the equivalence of (1) and (3) in~\cite[Theorem 3.1]{caranti}.
\end{proof}

\medskip

Recall that if $(B,+,\circ)$ is a skew brace, then $(B,+^{\operatorname{opp}},\circ)$ is a skew brace as well which goes under the name of {\it opposite} skew brace.

\begin{theo}\label{charQcirc}
Let $(B,+,\circ)$ be a non-trivial skew brace. If $(B,\circ)$ is torsion-free and locally cyclic, then $(B,+)=(C,+)\rtimes\langle a\rangle$, where $C$ is a trivial sub-brace with $(C,+)\simeq(B,\circ)$, and~$a$ is an additive involution inverting $(C,+)$. Moreover, $B$ is not Dedekind, $(C,+)$ is not $2$-divisible, and either
\begin{itemize}
    \item[\textnormal{(1)}] $\operatorname{Ker}(\lambda)=C$, $\lambda_a(b)=-b$ for every $b\in C$, and $\lambda_a(a)=a+a\circ a$, or
    \item[\textnormal{(2)}] $\operatorname{Ker}(\lambda)=B\ast C=\{0\}$, so in particular~\hbox{$a\circ b=a+b$} for every $b\in C$.
\end{itemize}

\noindent Conversely, let $(X,\circ)$ be any subgroup of $(\mathbb Q,+)$ which is not $2$-divisible, and choose $x\in X\setminus\{0\}$ with no square root in $X$.

\begin{itemize}

    \item[\textnormal{(1)}] Define a new operation $+$ on $X$ as follows: $x+x=0$, $b+d=b\circ d$, 
$x-d=x\circ d$, and $x+d+x=d^{-1}=-d$ for every $b,d\in\ X^2$. Then $(X,+,\circ)$ is a skew brace with $\operatorname{Ker}(\lambda)=X^2$.

    \item[\textnormal{(2)}] Define a new operation $+$ on $X$ as follows: $x+x=0$, $b+d=b\circ d$, 
$x+d=x\circ d$, and $x+d+x=d^{-1}=-d$ for every $b,d\in\ X^2$. Then $(X,+,\circ)$ is a skew brace with $\operatorname{Ker}(\lambda)=\{0\}$.

\end{itemize}

\noindent\textnormal(Note that the two skew braces are one the opposite of the other.\textnormal)
\end{theo}
\begin{proof}
Since $B$ is two-sided, $(B\ast B,+)$ is abelian (see Lemma 4.5 of \cite{nasy}) and (obviously) non-zero. Now, $B/B\ast B$ and~$B\ast B$ are radical rings, so their multiplicative and additive $0$-ranks coincide by~The\-o\-rem~B of~\cite{amberg}
--- the same result also shows that $(B,+)$ has finite abelian subgroup rank. Thus, $(B/B\ast B,+)$ is periodic while $(B\ast B,+)$ is torsion-free of rank $1$. 


If $(B\ast B,+)$ is central in $(B,+)$, then $(B,+)$ is nilpotent and the derived subgroup $(T,+)$ of~$(B,+)$ is periodic by Schur's theorem (see for example \cite{Ro72}, Corollary to Theorem 4.12). However, $(T,+)$ is covered by finite ideals (being of finite abelian subgroup rank and nilpotent), so $(T,\circ)$ is periodic and hence even zero. Thus,~$(B,+)$ is torsion-free abelian, and so also locally cyclic. By Theorem \ref{charQ+}, we obtain that $B$ is a trivial brace, a contradiction.

Assume $(B\ast B,+)$ is not central in $(B,+)$. Then its centralizer $(C,+)$ in $(B,+)$ (which is actually an ideal of $B$) has index $2$ in $(B,+)$. The case we dealt with in the previous paragraph yields that $C$ is a trivial brace. Since $(B,+)$ is not abelian, so $(B,+)=(C,+)\rtimes \langle a\rangle$, where $a$ has additive order $2$, the action of $a$ on $(C,+)$ is by inversion, and $(C,+)$ is a torsion-free locally cyclic group which is not $2$-divisible and isomorphic with $(B,\circ)$.

Now, if $\operatorname{Ker}(\lambda)\neq\{0\}$, then $(B,\circ)$ induces by $\lambda$-action a periodic group of automorphisms of $(B,+)$, which is only possible if this group of automorphisms has order $2$, so $\operatorname{Ker}(\lambda)=C$ because there is only one subgroup of index $2$ in $(B,\circ)$. Thus, either $\operatorname{Ker}(\lambda)=\{0\}$ or $\operatorname{Ker}(\lambda)=C$.


If $\operatorname{Ker}(\lambda)=C$, then $\lambda_a(b)=-b$ for every $b\in C$ and $\lambda_a(a)=a+c$ for some $c\in C$. But $\lambda_a(a)=-a+a\circ a$, and so $c=a\circ a$. Finally, assume $\operatorname{Ker}(\lambda)=\{0\}$,  and set $G=(B,+)\rtimes_\lambda(B,\circ)$. Clearly, $D=C_G(C,+)\cap (B,\circ)$ has index at most $2$ in~$(B,\circ)$, so either $D=(C,\circ)$ or $D=(B,\circ)$. However, if $D=(C,\circ)$, then $a^2\in\operatorname{Ker}(\lambda)=\{0\}$. Consequent\-ly,~$B\ast C=\{0\}$ and hence $a\circ b=a+\lambda_a(b)=a+b$ for every~\hbox{$b\in C$.}

The fact that $B$ is not Dedekind depends on the fact that $\langle a\rangle$ is a non-trivial sub-skew brace with a cyclic multiplicative group.

\medskip

Conversely, let $(X,+,\circ)$ be as in the statement. 

\smallskip

\noindent(1)\quad This follows at once from Lemma \ref{lemOld4.10}.


\bigskip


\noindent(2)\quad This is the opposite skew brace of that in (1). 
\end{proof}

\begin{cor}\label{thinfinitecyclic}
Let $(B,+,\circ)$ be a skew brace such that either $(B,+)$ or $(B,\circ)$ is torsion-free locally cyclic. Then the following conditions are equivalent\textnormal:
\begin{itemize}
    \item[\textnormal{(1)}] $B$ is trivial.
    \item[\textnormal{(2)}] $B$ is Dedekind.
    \item[\textnormal{(3)}] $(B,+)$ and $(B,\circ)$ are locally nilpotent.
\end{itemize}
\end{cor}
%
%


A further consequence of Theorem \ref{charQ+} is that if a skew brace has a torsion-free, locally cyclic, minimax additive group, then it must be trivial. Our next two results extend this observation.

\begin{theo}\label{divisibleideals}
Let $(B,+,\circ)$ be a skew brace. If either $(B,+)$ or $(B,\circ)$ is divisible periodic abelian of finite abelian subgroup rank, then $B$ is trivial.
\end{theo}
\begin{proof}
Suppose first that $(B,+)$ is divisible periodic abelian of finite abelian subgroup rank. Since $(B,+)$ is covered by finite characteristic subgroups, $(B,\circ)$ is periodic (actually, locally finite) and $B$ is index-preserving by Lemma \ref{analogcentralnilp}. We start by proving the statement when $(B,+)$ is a $p$-group for some prime $p$. In this case, $\big(B/\operatorname{Ker}(\lambda),\circ\big)$ is finite by Corollary to~Lem\-ma~3.28 of \cite{Ro72}. Then $\operatorname{Ker}(\lambda)$ has finite index also in~$(B,+)$ and so $B=\operatorname{Ker}(\lambda)$ because $(B,+)$ is divisible. The statement is thus proved in this case.

Now, we turn to the general case. For each prime $p$, let $B_p$ be the Sylow $p$-subgroup of~$(B,+)$. Then the previous paragraph shows that $B_p$ is a trivial sub-skew brace. Moreover, a further application of Corollary to Lemma 3.28 of \cite{Ro72} shows that $\big(B/\operatorname{Ker}(\lambda),\circ\big)$ is residually finite. It follows that each $B_p$ is contained in $\operatorname{Ker}(\lambda)$. But $\operatorname{Ker}(\lambda)$ is a sub-skew brace, and so~\hbox{$\operatorname{Ker}(\lambda)=B$.}

\medskip

Now, assume that $(B,\circ)$ is divisible periodic abelian of finite abelian subgroup rank. In this case, $B$ is two-sided. Suppose first that $(B,+)$ is abelian, so $B$ is a radical ring. Let $T$ be the torsion subgroup of $(B,+)$. Then~$T$ is an ideal of $B$, and Theorem B of \cite{amberg} yields that $(B/T,+)$ has $0$-rank $1$, while $(T,+)$ has finite abelian subgroup rank. It follows from~The\-o\-rem~\ref{charQ+} that $B=T$. Consequently, $(B,+)$ must be divisible and so the result follows from the first half of the proof.

Suppose $(B,+)$ is not abelian, so $B\ast B<B$ because $(B\ast B,+)$ is abelian by Lemma 4.5 of \cite{nasy}. Now, as in the previous paragraph, $(B\ast B,+)$ is periodic of finite abelian subgroup rank. Then $(B\ast B,+)$ is covered by finite characteristic subgroups, and the centralizers of these subgroups are characteristic subgroups of finite index of $(B,+)$. Thus, these centralizers coincide with $B$, and $B\ast B$ is central in $(B,+)$, which means that $(B,+)$ is nilpotent. Therefore, without loss of generality, we may assume that $(B,+)$ is a \v Cernikov $p$-group for some prime $p$. However, the finite residual of $(B,+)$ is an ideal of finite index, and consequently $(B,+)$ is abelian, a contradiction.
\end{proof}

\begin{theo}\label{thpolycyclic}
Let $(B,+,\circ)$ be a Dedekind skew brace.
\begin{itemize}
    \item[\textnormal{(1)}] If $(B,+)$ is soluble-by-finite and minimax, then $B$ is centrally nilpotent.
    \item[\textnormal{(2)}] If $(B,+)$ is torsion-free soluble-by-finite and minimax, then $B$ is a trivial brace.
    \item[\textnormal{(3)}] If $(B,\circ)$ is minimax abelian, then $B$ is centrally nilpotent. Moreover if $(B,\circ)$ is also torsion-free, then $B$ is a trivial brace.
\end{itemize}
\end{theo}
\begin{proof}
(1)\quad Let $(D,+)$ be the finite residual of $(B,+)$. Then $D$ is an ideal of $B$ which is trivial as a skew brace by Theorem \ref{divisibleideals}, and the quotient $(B/D,+)$ is residually finite.  If $(N,+)$ is any normal subgroup of~$(B,+)$ of finite index $n$, then $(nB,+)$ is a characteristic subgroup of~$(B,+)$ of finite index such that \hbox{$nB\subseteq N$.} Thus, $nB$ is an ideal of $B$ and $B/nB$ is centrally nilpotent by~Co\-rol\-la\-ry \ref{cordedcentralnilp}.

First, we show that $B/D$ is centrally nilpotent. To this aim, we may assume $D=\{0\}$. The structure of soluble-by-finite minimax groups and Theorem \ref{ThA} show that we may split the proof in two cases: $(B,+)$ abelian with a finite periodic subgroup; $(B,+)$ polycyclic-by-finite. However, in the latter case, every finite homomorphic image of $(B,+)$ is nilpotent, and hence~$(B,+)$ is nilpotent, so we only need to deal with the former.


Suppose $(B,+)$ is abelian. We may further assume $(B,+)$ is torsion-free because the periodic subgroup is finite. If $(B,+)$ is locally cyclic, then the statement follows from~The\-o\-rem~\ref{charQ+}. The\-o\-rem~10.34 of~\cite{Ro72} shows that $(B,+)$ is residually a finite $p$-group for infinitely many primes~$p$, so $B\ast B$ has infinite index in $B$, and induction on the Hirsch length yields that both $B\ast B$ and $B/B\ast B$ are centrally nilpotent. Then $B$ has a finite series of ideals whose factors are locally cyclic groups. Applying~Theorem \ref{ThA} shows that $B$ is centrally nilpotent.

Now, we show that $D\leq Z(B)$. Let $P$ be a Pr\"ufer $p$-group in $(D,+)$ for some prime~$p$, and write $D=P\times Q$ for some divisible subgroup $Q$ of $(D,+)$. If we are able to show that \hbox{$P\circ Q/Q\leq Z(B/Q)$,} then $P\leq Z(B)$ and  the arbitrariness of $P$ yields $P\leq Z(B)$. Thus, without loss of generality, we may assume that $(D,+)$ is a Pr\"ufer $p$-group for some prime $p$. Let $I/D$ be an ideal of $B/D$ such that $(I/D,+)$ is either infinite cyclic or cyclic of prime order. Being~De\-de\-kind,~$I/D$ is trivial, so $a\ast a\in D$ for a generator $a+D$ of $(I/D,+)$. Then factoring out $\langle a\ast a\rangle$, we may assume $a\ast a=0$, so $\langle a\rangle=\langle a\rangle_+$ is now a sub-skew brace and hence an ideal by hypothesis. An easy induction argument (on the lexicographically ordered set of pairs whose first component is the Hirsch length and the second is the order of the periodic subgroup of~$B/D$) shows that $\big(D+\langle a\rangle\big)/\langle a\rangle\leq Z\big(B/\langle a\rangle\big)$, so $D\leq Z(B)$ and the claim is proved.

\medskip

\noindent(2)\quad This follows from (1) and Corollary \ref{torsionfreeincredible}.


\medskip

\noindent(3)\quad In this case, $(B\ast B,+)$ is abelian (see~Lem\-ma~4.5 of~\cite{nasy}). Now, $B/B\ast B$ and $B\ast B$ are radical rings, so their additive groups are minimax abelian by~The\-o\-rem~A of \cite{amberg}. Then $(B,+)$ is metabelian and minimax, so $B$ is centrally nilpotent by (1). In the further assumption that $(B,\circ)$ is torsion-free, we have that the torsion subgroup of $(B,+)$ is zero, so $(B,+)$ is torsion-free and the statement follows from (2).
\end{proof}


We conclude this section with some consequences of the results established so far, concerning a natural nilpotency-like condition (which arises for example when $x\ast x=0$ for every element $x$ of the skew brace) that is often stronger than the Dedekind property.

\begin{theo}\label{generalizeadolf}
Let $(B,+,\circ)$ be a skew brace. Suppose that for each $b\in B$ which is an odd order or aperiodic additive/multiplicative element  either $\langle b\rangle_\circ$ or $\langle b\rangle_+$ is an ideal of~$B$.
Then $\langle b\rangle_+=\langle b\rangle_\circ$.

Moreover, the sub-skew brace $C$ generated by the elements of $B$ which have infinite additive \textnormal(or multiplicative\textnormal) order is a trivial brace.
\end{theo}
\begin{proof}
Let $b$ be an odd order or aperiodic additive/multiplicative element of $B$, and put \hbox{$H_b=\langle b\rangle_+$} or $H_b=\langle b\rangle_\circ$ if either $\langle b\rangle_+$ or $\langle b\rangle_\circ$ is an ideal of $B$. In any case, $H_b$ is Dedekind, and so centrally nilpotent by Theorem \ref{ThA}. If $H_b$ is infinite, then it is trivial by Corollary \ref{thinfinitecyclic}. If $H_b$ is finite of odd order, then the result follows at once from  Theorem \ref{cyclicideals}.

\medskip

Finally, let $a,b$ be additive aperiodic elements. Since $H_a$ is a trivial ideal of $B$, then $(pH_a,+)$ is an ideal of $B$ for each odd prime $p$. But then $a$ acts trivially on $H_b$ modulo $pH_a$ with respect to $+,\circ$ and $\lambda$. The arbitrariness of $p$ shows that $a$ acts trivially on $H_b$ with respect to $+,\circ$ and~$\lambda$. Consequently, $H_b$ is contained in the centre of the sub-skew brace $C$ of $B$ generated by the elements of $B$ which have infinite additive order. Therefore $C$ is a trivial brace.
\end{proof}

\medskip

The previous result fails for additive/multiplicative $2$-elements. In fact, there exist skew braces of order $4$ and $8$ with non-isomorphic additive and multiplicative groups but in which either all the additive or all the multiplicative subgroups generated by an element is always an ideal; for example, \texttt{SmallSkewbrace(8,5)} and \texttt{SmallSkewbrace(8,24)}. Thus, in order to obtain an analog of Theorem \ref{generalizeadolf} working for arbitrary skew braces, we need to add some hypothesis.


\begin{theo}\label{genadolfo2}
Let $(B,+,\circ)$ be a skew brace. Suppose $B$ satisfies the following properties:
\begin{itemize}
    \item If $b\in B$, then either $\langle b\rangle_\circ$ or $\langle b\rangle_+$ is an ideal of~$B$. 
    \item $(B,+)$ has a non-periodic element $x$.
\end{itemize} 

\noindent Then $B$ is a trivial brace.
\end{theo}
\begin{proof}
Clearly, the condition of Theorem \ref{generalizeadolf} is satisfied, so $\langle b\rangle_+=\langle b\rangle_\circ$ whenever $b$ is an odd order or aperiodic additive/multiplicative element of $B$. The same result shows that the sub-skew brace $C$ generated by the elements of $B$ having infinite additive order is trivial.

Let $b$ be an element of $B$ having finite additive and multiplicative order. Then $\langle b\rangle$ is a finite ideal of $B$ and such is also $\langle x\rangle$. Since $\langle b\rangle\cap\langle x\rangle=\{0\}$, $b+x$ is additively aperiodic. Consequently, $b\in C$ and hence $C=B$ is trivial.
\end{proof}

\begin{theo}\label{theoxastx}
Let $(B,+,\circ)$ be a Dedekind skew brace such that $(B,+)$ is not periodic. Let $X$ be a subset of $B$ such that $B=\langle X\rangle$. If $x\ast x=0$ for all $x\in X$, then $B$ is a trivial brace.
\end{theo}
\begin{proof}
If $x\in X$, then $\langle x\rangle=\langle x\rangle_+=\langle x\rangle_\circ$ is an ideal of $B$ (see Lemma \ref{bastb}). In particular, $X$ is both an additive and multiplicative set of generators for $B$. Thus, $B$ is hypercentral by~The\-o\-rem~\ref{ThA}. In order to prove the claim it is possible to assume $X$ is finite, so $B$ is actually centrally nilpotent. We prove the statement by induction on the central nilpotency class of $B$. If $B=Z(B)$, then the statement is obvious. Assume $B=Z_n(B)$ for some $n$. Then~$B/Z(B)$ is a trivial brace and hence $B=Z_2(B)$.

Since the elements of finite additive order in $B$ form a finite subgroup, there must be $y\in X$ of infinite additive order. Moreover, $\langle y\rangle/\langle y^p\rangle\leq Z\big(B/\langle y^p\rangle\big)$ for every prime $p$, so $y\in Z(B)$. Now, if~$x$ is any element of $X$, then $$
\begin{array}{c}
(x+y)\ast (x+y)=(x+y)\ast x+(x+y)\ast y=(x\circ y)\ast x+(x\circ y)\ast y\\[0.2cm]
=x\ast (y\ast x)+x\ast x+y\ast x+x\ast (y\ast y)+x\ast y+y\ast y=0.
\end{array}
$$ Consequently, if $x$ has finite additive order, then $x+y$ has infinite additive order and the ideal $\langle x+y\rangle=\langle x+y\rangle_+=\langle x+y\rangle_\circ$ is contained in $Z(B)$, which means that $x\in Z(B)$. In any case, $X$ lies in $Z(B)$, and hence $B=Z(B)$ is a trivial brace.
\end{proof}

\begin{cor}
Let $(B,+,\circ)$ be a Dedekind skew brace such that $(B,+)$ is not periodic. If $b\ast b=0$ for all $b$ in $B$, then $B$ is a trivial brace.
\end{cor}

\begin{cor}\label{corsolution}
Let $(X,r)$ be a solution of the YBE such that $G(X,r)$ is a Dedekind skew brace. If~\hbox{$r(x,x)=(x,x)$} for every $x\in X$, then $(X,r)$ is a twist solution. In particular, $G(X,r)$ is a group which is isomorphic to~$\mathbb Z^X$.
\end{cor}
\begin{proof}
The commutativity of the diagram \eqref{diagram3} and the hypothesis yield that \hbox{$\iota_X(x)\ast \iota_X(x)=0$} for every $x\in X$. Consequently, Theorem \ref{theoxastx} shows that $G(X,r)$ is trivial, which means that $(X,r)$ is the twist solution.
\end{proof}

\vspace{3cm}
\noindent {\bf Acknowledgments}\quad  Ferrara and Trombetti are members of the non-profit association ‘‘AGTA --- Advances in Group Theory and Applications’’ (www.advgrouptheory.com). Funded by the European Union - Next Generation EU, Missione 4 Componente 1 CUP B53D23009410006, PRIN 2022 --- 2022PSTWLB --- Group Theory and Applications. 
Del Corso gratefully acknowledges the support of the MIUR Excellence Department Project awarded to the Department of Mathematics, University of Pisa, CUP I57G22000700001.
Del Corso has performed this activity in the framework of the PRIN 2022, title ‘‘Semiabelian varieties, Galois representations and related Diophantine problems’’. All the authors are supported by \hbox{GNSAGA} (INdAM).


\newpage

\begin{flushleft}
\rule{8cm}{0.4pt}\\
\end{flushleft}

\bigskip
\bigskip

{
\sloppy
\noindent
Andrea Caranti

\noindent
Dipartimento di Matematica

\noindent
Università degli Studi di
Trento

\noindent
via Sommarive 14, I-38123 Trento (Italy)

\noindent
e-mail: andrea.caranti@unitn.it

\noindent
url: https://caranti.maths.unitn.it/
}

\bigskip
\bigskip

{
\sloppy
\noindent
Ilaria Del Corso

\noindent
Dipartimento di Matematica

\noindent
Università di Pisa

\noindent
Largo Bruno Pontecorvo, 5, 56127 Pisa (Italy)

\noindent
e-mail: ilaria.delcorso@unipi.it

}

\bigskip
\bigskip

{
\sloppy
\noindent
Massimiliano di Matteo

\noindent
Dipartimento di Matematica e Fisica

\noindent
Università degli Studi della Campania  ``Luigi Vanvitelli''

\noindent
viale Lincoln 5, Caserta (Italy)

\noindent
e-mail: massimiliano.dimatteo@unicampania.it 
}

\bigskip
\bigskip

{
\sloppy
\noindent
Maria Ferrara

\noindent 
Dipartimento di Matematica e Fisica

\noindent
Università degli Studi della Campania “Luigi Vanvitelli”

\noindent
viale Lincoln 5, Caserta (Italy)






\noindent
e-mail: maria.ferrara1@unicampania.it 

}

\bigskip
\bigskip

{
\sloppy
\noindent
Marco Trombetti

\noindent 
Dipartimento di Matematica e Applicazioni ``Renato Caccioppoli''

\noindent
Università di Napoli Federico II

\noindent
Complesso Universitario Monte S. Angelo

\noindent
Via Cintia, Napoli (Italy)

\noindent
e-mail: marco.trombetti@unina.it 

}



\newpage
\begin{thebibliography}{10}


\bibitem{acri2020skew}
E. Acri and M. Bonatto:
‘‘Skew braces of size $pq$’’,
{\it Comm. Algebra} 48 (2020), 1872--1881.



\bibitem{amberg}
B. Amberg and O. Dickenschied: ‘‘On the adjoint group of a radical ring’’,
{\it Canad. Math. Bull.} 38 (1995), no. 3, 262--270.


\bibitem{p3}
D. Bachiller: ‘‘Classification of braces of order $p^3$’’, {\it J. Pure Appl. Algebra} 219 (2015), no. 8, 3568--3603.


\bibitem{centnilpo} 
A. Ballester-Bolinches, R. Esteban-Romero, M. Ferrara, V. Pérez-Calabuig and M. Trombetti: ‘‘Central nilpotency of left
skew braces and solutions of the Yang–Baxter equation’’,
{\it Pacific J. Math.} 335 (2025), No. 1, 1--32.


\bibitem{supersol} 
A. Ballester-Bolinches, R. Esteban-Romero, M. Ferrara, V. Pérez-Calabuig and M. Trombetti: ‘‘Finite skew braces of square-free order and supersolubility’’, {\it Forum Math. Sigma} 12 (2024), e39.

\bibitem{poly} 
A. Ballester-Bolinches, R. Esteban-Romero, M. Ferrara, V. Pérez-Calabuig and M. Trombetti: ‘‘The Word Problem and the Conjugacy Problem for the structure skew brace of a solution of the Yang--Baxter Equation are solvable’’, to appear.


\bibitem{ballester2024soluble} 
A. Ballester-Bolinches, R. Esteban-Romero, L.A. Kurdachenko, V. Pérez-Calabuig: ‘‘Soluble skew left braces and soluble solutions of the Yang--Baxter Equation’’, {\it Adv. Math.} 455 (2024), 109880.

\bibitem{ballesterDedekindii} 
A. Ballester-Bolinches, R. Esteban-Romero, L.A. Kurdachenko, V. Pérez-Calabuig: ‘‘On left braces in which every subbrace is an ideal II’’, to appear.


\bibitem{dedekindbraces} 
A. Ballester-Bolinches, R. Esteban-Romero, L.A. Kurdachenko, V. Pérez-Calabuig: ‘‘On left braces in which every subbrace is an ideal’’, {\it Results Math.} 80 (2025), Article Number 21, 21pp.


\bibitem{edinburgh} 
A. Ballester-Bolinches, M. Ferrara, V. Pérez-Calabuig and M. Trombetti: ‘‘A note on right-nil and strong-nil skew braces’’,
{\it Proc. A Royal Soc. Edinburgh}; 
doi:10.1017/prm.2024.133.

\bibitem{Baxter}
R.J. Baxter: ‘‘Partition function of the eight-vertex lattice model’’, {\it Ann. Physics} 70 (1972), 193--228. 


\bibitem{byott} N. Byott: ‘‘Hopf-Galois structures on almost cyclic field extensions of $2$-power degree’’, {\it J. Algebra} 318 (2007), no. 1, 351--371.

\bibitem{campedel} E. Campedel, A. Caranti, and I. Del Corso: 
‘‘Hopf-Galois structures on extensions of degree $p^2q$ and skew braces of order $p^2q$: the cyclic Sylow $p$-subgroup case’’, {\it J. Algebra} 556 (2020), 1165--1210.

\bibitem{caranti} A. Caranti: ‘‘Bi-skew braces and regular subgroups of the holomorph’’, {\it J. Algebra} 562 (2020), 647--665.


\bibitem{stefanelli}
F. Catino, I. Colazzo, and P. Stefanelli:
\emph{Skew left braces with non-trivial annihilator},
{\it J. Algebra Appl.} 18 (2019), no. 2, 1950033, 23 pp.

\bibitem{cedo1}
F.~Cedó, E.~Jespers, Ł.~Kubat, A.~Van Antwerpen, and C.~Verwimp: ‘‘On various types of nilpotency of the structure monoid and group of a set-theoretic solution of the Yang--Baxter equation’’ {\em J. Pure Appl. Algebra} 227 (2023), 107194.

\bibitem{nilpotent}
F. Ced\'{o}, A. Smoktunowicz, and L. Vendramin: ‘‘Skew left braces of nilpotent type’’, {\it Proc. Lond. Math. Soc.} (3) 118 (2019), no. 6, 1367--1392.

\bibitem{ChildsYBE}
L.N. Childs: ‘‘Skew left braces and the Yang-Baxter equation’’, {\it New York J. Math.} 30 (2024), 649--655.

\bibitem{colazzo}
I. Colazzo, M. Ferrara, and M. Trombetti: ‘‘On derived-indecomposable solutions of the Yang--Baxter equation’’, {\it Publ. Math.} 69 (2025), no. 1, 171--193.

\bibitem{idc2023}
I. Del Corso: ‘‘Module braces: relations between the additive and the multiplicative groups’’,
{\it Ann. Mat. Pura Appl.} 202  (2023), 3005--3025.

\bibitem{GAP}
The GAP Group: ‘‘GAP – Groups, Algorithms, and Programming’’, Version 4.12.2 (2022).

\bibitem{GV}
L. Guarnieri and L. Vendramin: ‘‘Skew braces and the Yang--Baxter equation’’, {\it Math. Comp.} 86 (2017), no. 307, 2519--2534.

\bibitem{JVV}
E. Jespers, A. Van Antwerpen, and L. Vendramin: ‘‘Nilpotency of skew braces and multipermutation solutions of the Yang-Baxter equation’’, {\it Commun. Contemp. Math.} 25 (2023), no. 9, Paper No. 2250064, 20 pp.


\bibitem{Kohl} T. Kohl: ‘‘Classification of the Hopf Galois structures on prime power radical extensions’’, {\it J. Algebra} 207 (1998), no. 2, 525--546.

\bibitem{nasy} T. Nasybullov: ‘‘Connections between properties of the additive and the multiplicative groups of a two-sided skew brace’’, {\it J. Algebra} 540 (2019), 156--167.

\bibitem{pompili} M. Pompili: ‘‘Yang-Baxter Equation
and Category of Skew Braces’’, {\it M.Sc. Thesis}, Padova (2022).

\bibitem{Ro72} D.J.S. Robinson: ‘‘Finiteness Conditions and Generalized Soluble Groups’’, {\it Springer}, Berlin (1972).


\bibitem{acourse}
D.J.S. Robinson: ‘‘A Course in the Theory of Groups’’, 2ed., {\it Springer}, New York (1996). 

\bibitem{Rump0}
W. Rump: ‘‘Braces, radical rings, and the quantum Yang-Baxter equation’’, {\it J. Algebra} 307 (2007), no. 1, 153--170.

\bibitem{rump2007classification}
W. Rump: ‘‘Classification of cyclic braces’’, {\it J. Pure  Appl. Algebra} 209 (2007), no. 3, 671--685.

\bibitem{SV}
A. Smoktunowicz and L. Vendramin: ‘‘On skew braces (with an appendix by N. Byott and L. Vendramin)’’, {\it J. Comb. Algebra} 2 (2018), no. 1, 47--86.

\bibitem{stefanellotrap} L. Stefanello and S. Trappeniers: ‘‘On bi-skew braces and brace blocks’’, {\it J. Pure Appl. Algebra} 227 (2023), 107295.

\bibitem{stefanello-tesi} L. Stefanello: ``Hopf--Galois Structures, Skew Braces, and Their Connection'', PhD. Thesis, {\it University of Pisa}, Pisa (2024).


\bibitem{trombetti2023structure} M. Trombetti: ‘‘The structure skew brace associated with a finite non-degenerate solution of the Yang–Baxter equation is finitely presented’’, {\it Proc. AMS} 152(2) (2024), 573--583


\bibitem{cindy} C. Tsang: ‘‘Analogs of the lower and upper central series in skew braces: a survey’’, {\it Commun. Math.} 33 (2025), no. 3, Paper No. 11, 30 pp.


\bibitem{Yang}
C.N. Yang: ‘‘Some exact results for the many-body problem in one dimension with repulsive delta-function interaction’’, {\it Phys. Rev. Lett.} 19 (1967), 1312--1315.








\end{thebibliography}

\end{document}